\documentclass{amsart}

\usepackage{latexsym}
\usepackage{url}

\usepackage{amsmath,graphicx,bbm}
\usepackage{amsthm,amssymb,verbatim}
\usepackage{mathrsfs}
\usepackage[footnotesize,bf]{caption}
\usepackage[left=1.25in,right=1.25in,top=1in]{geometry}
\usepackage{tikz,pgfplots,pgfplotstable}
\usepackage[textwidth=1.25in,textsize=tiny]{todonotes}
\usepackage{setspace}
\usepackage{algorithmic} % For algorithm listings
\usepackage{algorithm} % For algorithm listings
\usepackage{color}
\usepackage{subfigure}

\usepackage{acncommands}
\usepackage{tikz,pgfplots,pgfplotstable}
\usepgfplotslibrary{external}
\begin{document}

\newtheorem{thm}{Theorem}
\parskip=.1in

\title[Approximating the WP Geodesics in $T(1)$ by Teichons]{Approximating the Weil-Petersson Metric Geodesics on the Universal Teichm\"uller space by Singular Solutions}         % Enter your title between curly braces
\author{Sergey Kushnarev}
\thanks{Dr Sergey Kushnarev, Research Fellow at Computational Functional Anatomy Laboratory, Department of Bioengineering, National University of Singapore, 2 Engineering Drive 3, Block E3 \#05-30, Singapore. Email: bieserge@nus.edu.sg}
\author{Akil Narayan}
\thanks{Akil Narayan, Assistant Professor, Mathematics Department, University of Massachusetts Dartmouth, 285 Old Westport Road, Dartmouth, MA, USA. Email: akil.narayan@umassd.edu.}        % Enter your name between curly braces
\date{January 2012}          % Enter your date or \today between curly braces
\begin{abstract}

We propose and investigate a numerical shooting method for computing geodesics in the Weil-Petersson ($WP$) metric on the universal Teichm\"uller space $T(1)$. This space, or rather the coset subspace $\PSL_2(\R)\backslash\Diff(S^1)$, has another realization as the space of smooth, simple closed planar curves modulo translations and scalings. This alternate identification of $T(1)$ is a convenient metrization of the space of shapes and provides an immediate application for our algorithm in computer vision. The geodesic equation on $T(1)$ with the $WP$ metric is EPDiff($S^1$), the Euler-Poincare equation on the group of diffeomorphisms of the circle $S^1$, and admits a class of soliton-like solutions named Teichons \cite{kushnarev_teichons:_2009}. Our method relies on approximating the geodesic with these teichon solutions, which have momenta given by a finite linear combination of delta functions. The geodesic equation for this simpler set of solutions is more tractable from the numerical point of view. With a robust numerical integration of this equation, we formulate a shooting method utilizing a cross-ratio matching term. Several examples of geodesics in the space of shapes are demonstrated.
\end{abstract}
\maketitle

%%edited oct-19-2009
\def\R{\mathbb{R}} % real numbers
\def\C{\mathbb{C}} % complex numbers
\def\Z{\mathbb{Z}} % integers
\def\D{\mathbb{D}} % unit disk
\def\Diff{\mathbf{Diff}} % group of diffeomorphisms
\def\Homeo{\mathbf{Homeo}} % group of homeomorphisms
\def\Vec{\mathbf{Vec}} % algebra of smooth vector fields
\def\PSL{\mathrm{PSL}} % projective linear group
\def\Ad{\mathrm{Ad}} % Adjoint action
\def\ad{\mathrm{ad}} % adjoint action
\def\g{\mathfrak{g}} % Lie algebra g
\def\h{\mathfrak{h}} % Lie algebra h
\def\p{\mathfrak{p}} % Lie algebra p
\def\bs{\backslash} % backslash
\def\cnst{\mathrm{c}} % math constant in formulas
\def\H{\mathcal{H}} % Hilbert transform
\def\th{\theta}
\def\dth{\Delta\theta} % Delta theta
\newcommand{\abs}[1]{{|{#1}^3-{#1}|}} % |[1]^3-[1]|
\newcommand{\adtr}[2]{{\ad^\dag_{#1}#2}} % transpose of adjoint
\def\Hor{Hor} % Horizontal bundle
\def\const{\text{const}} % constant
\newcommand{\bpm}{\begin{pmatrix}}
\newcommand{\epm}{\end{pmatrix}}

%\tableofcontents

\section{Introduction}      \label{s:intro}
 We consider the Weil-Petersson (WP) metric on the coset space $\PSL_2(\R)\backslash\Diff(S^1)$. This coset space (or its completion in the WP metric or in the Teichm\"uller topology) is known as the universal Teichm\"uller space and is well-known in many contexts: in the classification of Riemann surfaces \cite{hubbard_teichmuller_2006}, conformal and quasi-conformal maps \cite{lehto_univalent_1986}, string theory \cite{bowick_string_1987} and most recently computer vision \cite{sharon_2d-shape_2006}. Its completion in the WP metric is an infinite dimensional homogeneous complex K\"ahler-Hilbert manifold \cite{teo_weil-petersson_2006}. 

 As we will explain in Section \ref{s:T1} below, a particular dense subset of the universal Teichm\"uller space $T(1)$ is given by $\PSL_2(\R)\bs\Diff(S^1)$, where $\Diff(S^1)$ is the group of $C^\infty$ diffeomorphisms of $S^1$, and $\PSL_2(\R)$ is a subgroup of the M\"obius selfmaps of the unit disk, see \eqref{eq:psl2-definition} and the surrounding discussion. This coset space is a Riemannian manifold for the WP metric and has another realization as the space of smooth simple closed curves modulo translations and scalings. (Therefore we will hereafter use the terms `shape', `diffeomorphism', `fingerprint', or `welding map' to refer to members in this dense subset of $T(1)$.) Endowing a shape space with a metric and computing geodesics between shapes is a central problem in computer vision. It aids in recognition and classification, enables computational strategies to address the clique problem, and allows us to perform statistics on shapes. Another application is in the emerging field of computational anatomy, where quantitative analysis of anatomical variability is important \cite{grenander_computational_anatomy}. It is this major application that provides the motivation for this work. 

There are several advantages to using the Weil-Petersson metric to compare two-dimensional shapes. First, any two smooth shapes can be connected with a Weil-Petersson geodesic \cite{gay-balmaz_geometry_2009}. Second, all sectional curvatures of the metric are negative \cite{teo_weil-petersson_2006}. Thus geodesics connecting two shapes are unique \cite{karcher}. We are not aware of any other metric currently used in the pattern theory literature that has these properties. Ultimately, the uniqueness of geodesics allows one to perform consistent statistical analysis on shape databases via the initial momentum representation of the shape, but this application is beyond the scope of this paper.

%The main contribution of this paper is a numerical algorithm for computation of geodesics on this space. Our method finds a geodesic connecting two endpoints in this metric by approximating it with a finite linear combination of soliton-like solutions. The soliton-like properties of these solutions allow more efficient numerical computation of geodesics. With a robust method for forward propagation on geodesics, we formulate a boundary value shooting problem, requiring an appropriate matching term.

%Therefore the solutions of this equation can also be thought of as paths in the space of simple closed plane curves which minimize a certain energy. The soliton property means that, in a certain sense, their momentum is concentrated at a finite set of points.

Geodesic equations of groups of diffeomorphisms on a manifold $M$ were first studied in Arnold's ground-breaking paper \cite{arnold_sur_1966}. Arnold considered in particular the group of volume preserving diffeomorphisms of Euclidean space in its $L^2$ metric and found the geodesic equation for the vector field $\vec v(\vec x,t)$ to be Euler's fluid flow equation (see \cite{arnold_topological_1998} for a full exposition). Other examples include the periodic Korteweg-deVries (KdV) equation and the periodic Camassa-Holm (C-H) equation \cite{camassa_integrable_1993}. These equations are geodesic equations on the Virasoro group, a central extension by $S^1$ of the group $\Diff(S^1)$ of the diffeomorphisms of $S^1$, for the $L^2$ and $H^1$ metric respectively. KdV and C-H are two completely integrable partial differential equations and have soliton solutions. 
Holm and collaborators have found that the geodesic equation on $\Diff(\mathbf R^n)$ admits special solutions with many of the properties of solitons: for each fixed time, they are diffeomorphisms which are largely localized in space and retain their general shape as they evolve; furthermore they interact somewhat like KdV solitons \cite{holm_soliton_2001}. There are not, however, infinitely many conserved quantities so they are not true solitons.

Singular solutions first arose as peakons (from `peaked solitons') for a completely integrable Hamiltonian water wave equation, C-H in \cite{camassa_integrable_1993}. The peaks occurred where the velocity profiles of the C-H equation had discontinuity in its slope. These peaks correspond to Dirac delta distributions of the associated momentum. The EPDiff equation for other metrics was later found independently in \cite{trouve_1995}, and its singular solutions were shown to be important as landmarks in shape analysis  \cite{grenander_computational_anatomy, mumford_book}. Later they were shown to comprise a singular momentum map for the right action of the diffeomorphisms on embeddings in any dimension \cite{holm_momentum_2005}. Currently, the use of EPDiff and its landmark solutions is standard in shape analysis \cite{holm_etal_2004, miller_etal_2002, miller_geodesic_2006}.

It turns out that considering the Weil-Petersson metric on the coset space $\PSL_2(\R)\bs\Diff(S^1)$ yields another example of a geodesic equation that is similar to KdV and C-H. This equation describing evolution of the velocity field $v(t,\theta)$ is
\begin{equation}\label{eq:epdiff}
m_t + 2m v_\th + v m_\th = 0,\text{ where } m = -\H(v_\th+v_{\th \th \th}),
\end{equation}
and $\H$ is the periodic Hilbert transform defined by convolution with $\tfrac{1}{2\pi}\text{ctn}(\th/2)$. 
%Here, we may invert the relationship between $m$ and $u$ and write:
%$$ 
%v(\theta,t) = \int_{S^1} G(\theta-\xi)m(\xi,t)d\xi = G*m. 
%$$
%The integral kernel, or Green's function, $G(\theta)$ turns out to be given in the Fourier domain by 
%$$ 
%G(\th) = 2\sum_{n=2}^\infty \frac{\cos(n\th)}{(n^3-n)}.
%$$
%Note that $m(\cdot,t)$ is always orthogonal to 1, cos and sin.

It is not known if \eqref{eq:epdiff} is completely integrable but it admits a class of soliton-like solutions which we consider in this paper: solutions in which $m$ can be represented as a finite sum of weighted Dirac delta functions. Darryl Holm suggested the portmanteau {\em teichons} to describe these soliton-like solutions on Teichm\"uller space and their corresponding geodesics. We adopt this terminology in this paper.

We use an $N$-teichon ansatz (a sum of $N$ teichons) to reduce the integro-differential equation \eqref{eq:epdiff} to a finite-dimensional system of ordinary differential equations. In this way we approximate geodesics between any two points in $\PSL_2(\R)\bs\Diff(S^1)$ with $N$-teichon geodesic evolutions. We use this teichon formulation to shoot from an initial shape to a terminal shape that must then be compared with the target shape. Because we are considering the coset space $\PSL_2(\R)\bs\Diff(S^1)$, using a standard matching term on diffeomorphisms is not possible: the quotient space ambiguity prevents straightforward comparisons of diffeomorphisms (e.g. pointwise matching). We address this difficulty by using cross-ratios in the matching term; their invariance within an equivalence class on $\PSL_2(\R)\bs\Diff(S^1)$ allows for accurate matching on the coset space.

One existing geodesic computing algorithm, described in \cite{sharon_2d-shape_2006} and investigated in \cite{kushnarev_geometry_2010}, suffers from several limitations. In particular there are numerical difficulties in matching shapes that are not close to a circular shape. The current approach aims to mitigate this shortcoming. A recent approach given in \cite{feiszli_2012} is more competitive with our approach.

This paper is organized as follows. Section \ref{s:T1} introduces the background on universal Teichm\"uller space $T(1)$, fingerprints (also called welding maps), and the Weil-Petersson metric. With the $WP$ metric, Section \ref{s:epdiff} discusses the geodesic equation on the Teichm\"uller space (also known as EPDiff), and Section \ref{s:teichons} discusses teichon solutions of EPDiff. In Section \ref{s:nummethod} we describe the shooting method and the matching functional used for shooting, along with details of the gradient computation. In Section \ref{s:numres} we demonstrate the utility of our method with several examples.

\section{Shapes as diffeomorphisms of the circle $S^1$}  \label{s:T1}
\subsection{Fingerprints}\label{ss:modA}
Let $\D_{\mathrm{int}}$ be the open unit disk in the complex plane $\C$, i.e.\ 
$\D_{\mathrm{int}} = \{z\in\C\mid |z|<1\}$, and let $\D_{\mathrm{ext}} = \{z\in\C\mid |z|>1\}$ be its exterior. 
For every simple closed curve $\Gamma$ in $\C$ denote by $\Gamma_{\mathrm{int}}$ its union with the region enclosed by it, and denote by $\Gamma_{\mathrm{ext}}$ its union with the infinite region outside of $\Gamma$ (including $\infty$).

Then by the Riemann mapping theorem, for all $\Gamma$ there exist two conformal maps
\begin{align*}
f_\mathrm{int}&: \D_\mathrm{int} \to \Gamma_\mathrm{int},\\
f_\mathrm{ext}&: \D_\mathrm{ext} \to \Gamma_\mathrm{ext}.
\end{align*}
The interior map $f_\mathrm{int}$ is unique up to replacing $f_\mathrm{int}$ by $f_\mathrm{int}\circ A$ for any M\"obius transformation $A:\D_\mathrm{int} \to \D_\mathrm{int}$, where $A$ defined as 
\begin{align}\label{eq:psl2-definition}
A(z) = \frac{az+b}{\bar{b}z+\bar{a}}, \ |a|^2-|b|^2 = 1.
\end{align}
This subgroup of M\"obius group of selfmaps of the circle is denoted $\PSL_2(\R)$.

The map $f_\mathrm{ext}$ is chosen uniquely via the following normalization: we choose a unique M\"obius map $A$, such that $f_\mathrm{ext}\circ A$ maps $\infty$ to $\infty$, and that its differential carries the real positive axis of the $\D$-plane at infinity to the real positive axis of the $\Gamma$-plane at infinity. Thus the ambiguity in the choice of $f_\mathrm{ext}$ is eliminated for every $\Gamma$.

The goal of this construction is to define the map $\psi$ which is called the `fingerprint' (in the Teichm\"uller theory this is known as a `welding map') of the shape
\begin{align}\label{eq:weld-definition}
\psi = f_\mathrm{int} ^{-1} \circ f_\mathrm{ext} \in \PSL_2(\R)\backslash \Diff(S^1).
\end{align}
Note, that $f_\mathrm{ext}(S^1)=\Gamma,\ f_\mathrm{int}^{-1}(\Gamma)=S^1$. The fingerprint $\psi: S^1 \to S^1$ is a real-valued orientation-preserving diffeomorphism, and it uniquely identifies the shape $\Gamma$ (modulo scaling and rigid translations). Due to the M\"obius transformation ambiguity in the choice of $f_\mathrm{int}$, we see by construction that $\psi$ is a member of the right coset space $\PSL_2(\R)\backslash \Diff(S^1)$.
%, where $\PSL_2(\R)$ is a subgroup of M\"obius maps that preserves the unit circle.
An example of a shape along with four realizations of its fingerprint is given in Figure \ref{fig:welding-maps}.

The inverse map from diffeomorphisms to shapes is defined as follows: starting with $\psi$, construct an abstract Riemann surface by `welding' the boundaries of $\D_\mathrm{int}$ and $\D_\mathrm{ext}$ via $\psi$. The resulting Riemann surface must be conformally equivalent to the Riemann sphere. Choose a conformal map $f$ from the welded surface to the sphere taking $\infty \in \D_\mathrm{ext}$ to itself and having real positive derivative there. Let $\Gamma = f(S^1)$ (for details and the numerical implementation see \cite{sharon_2d-shape_2006}). 

One can equally well define the fingerprint to be 
$$ 
\psi^i=f_\mathrm{ext} ^{-1} \circ f_\mathrm{int} \in \Diff(S^1) / \PSL_2(\R),
$$
which is simply the inverse of our fingerprint. This alternate version is the definition used in \cite{sharon_2d-shape_2006}. However, in this paper we choose {\it right} cosets and put the M\"obius ambiguity on the left.

\subsection{Weil-Petersson Norm on the Lie Algebra of $\Diff(S^1)$}
The Lie algebra of the group $\Diff(S^1)$ is given by the vector space $\Vec(S^1)$ of smooth periodic vector fields $v(\theta)\partial/\partial\theta$ on the circle. In \cite{nag_diff_1990} it has been shown that the embedding $\PSL_2(\R)\bs\Diff(S^1)\hookrightarrow T(1)$ is holomorphic and the pullback of the Weil-Petersson metric for $v \in \Vec(S^1)$ can be expressed as:
\begin{align*}
\|v\|_{WP}^2 &= \sum_{n\in \widehat\Z}|n^3-n| |v_n|^2\\
	     &= \int_{S^1} Lv(\th) v(\th) d\th.
%\label{WPmetric}
\end{align*}
Here $v(\theta)=\sum_{n=-\infty}^{\infty}v_n e^{in\theta}$ (where $\overline{v_n}=v_{-n}$ for the vector field to be real), and $\widehat\Z = \Z\backslash\{n=0,\pm 1\}$. The Weil-Petersson operator $L$ is an integro-differential operator and it has the form
\begin{align}\label{eq:wp-operator}
L = -\H(\partial_\th^3+\partial_\th).
\end{align}
Above, $\H$ is the periodic Hilbert transform, defined as a convolution with $\frac{1}{2\pi}\cot(\th/2)$.

The null space of the $L$ operator is given by the vector fields whose only Fourier coefficients are $v_{-1}, v_0$ and $v_1$, i.e. vector fields of the type $(a+b\cos\theta+c\sin\theta)\partial/\partial\theta$. These vector fields are exactly in the Lie algebra $\mathfrak{sl}_2(\R)$ of the Lie group $\PSL_2(\R)$.
%$psl_2(\R)$ of the Lie group $\PSL_2(\R)$.

%$$
%\|\Ad_A(v)\|_{WP} = \|v\|_{WP}.
%$$

\subsection{Extending the WP Metric to $\PSL_2(\R)\bs\Diff(S^1)$}
Consider any Lie group $G$, a subgroup $H$, and let $\g$ and $\h$ be their corresponding Lie algebras. 

Any norm $\|\cdot\|$ on the Lie algebra of $G$ which is zero on the Lie subalgebra of $H$ and which satisfies $\|\Ad_h(v)\|=\|v\|$ for all $h\in H$ induces a Riemannian metric on coset space $H\backslash G$ which is invariant by all right multiplication maps $R_g: H\backslash G \to H\backslash G$, $g\in G$ \cite{kushnarev_geometry_2010}.

In particular this applies to $G = \Diff(S^1), H=\PSL_2(\R)$ and the above WP norm on vector fields, hence it gives the right-invariant {\it WP-Riemannian metric} on the coset space $\PSL_2(\R)\bs\Diff(S^1)$. 

Consider any two diffeomorphisms $\psi_0, \psi_1 \in \uts$. The Riemannian distance induced by the WP norm on vector fields is given by 
\begin{align}\label{eq:riemannian-distance}
  L = \int_0^1 \| v(s)\|_{W P} \dx{s},
\end{align}
where $v$ is a vector field that carries $\psi_0$ to $\psi_1$:
\begin{align}\label{eq:v-evolution}
  \phi_t(\theta) = \psi_0(\theta) + \int_0^t v(\theta, s) \dx{s}
\end{align}
This above notation extends for the remainder of this paper: $\psi_1(\theta)$ and $\psi_0(\theta)$ are the initial and target welds (shapes), and the path $\phi_t(\theta) = \phi(\theta,t)$ is a geodesic flow. 
Vector fields $v$ that minimize the distance \eqref{eq:riemannian-distance} are geodesics on $\uts$, and it is a standard fact from variational calculus that vector fields $v$ corresponding to geodesics satisfy $\|v(s,\cdot)\|_{W P} = \mathrm{const}$.

\section{The geodesic equation}\label{s:epdiff}

The Euler-Poincar\'e equation for diffeomorphisms (hereafter `EPDiff') 
is a variant of Euler's equations for fluid flow. It describes geodesics on the Lie group of diffeomorphisms of $\R^n$ in any right invariant metric given on vector fields by $\|v\|^2 = \int_{\R^n} \langle Lv,v\rangle dx$ for some positive definite self-adjoint operator $L$ (where $\langle \cdot, \cdot \rangle$ is the canonical $L^2$ pairing). The general EPDiff($\R^n$) is derived in \cite{arnold_sur_1966} and has the form
$$
\frac{\partial}{\partial t}Lv + (v\cdot\nabla)(Lv) + \mathrm{div} v\; Lv + Dv^t\cdot Lv = 0,
$$
where $v$ is a smooth vector field in $\R^n$, $\nabla=\left(\frac{\partial}{\partial x_1},\ldots,\frac{\partial}{\partial x_n}\right)^T$ is the divergence operator, $L$ is a self-adjoint differential operator and $Dv$ is a Jacobian matrix.
%\red{Sergey: what is all this notation? $D v = ?$, $v^t = ?$, $v.w = ?$ $v\cdot w = ?$, $\nabla$ vs $D$?}

The space $\uts$ that interests us is not a group and is instead a homogeneous space, 
%In our case, we have the homogeneous space $\PSL_2(\R)\bs\Diff(S^1)$, not a group. 
but it has been shown in \cite{khesin_euler_2003,kushnarev_geometry_2010} that Arnold's formula for geodesics on Lie groups extends to the case of a homogeneous spaces $H\bs G$.

 Given a path $\phi_t(\theta) = \phi(\th,t)$ in $\Diff(S^1)$, let $v(\theta,t)=\frac{\partial\phi}{\partial t}(\phi^{-1}(\theta,t),t)$ be the scalar vector field it defines on a circle and let $L$ be the Weil-Petersson differential operator $L=-\H(\partial_\theta^3+\partial_\theta)$. Then EPDiff takes the form
\begin{equation}
(Lv)_t + v.(Lv)_\theta + 2v_\theta.Lv = 0.
\label{EPDiff}
\end{equation}
Above, $v(\theta,t)$ is called the velocity of the path, $m(\theta,t) = Lv(\th,t)$ is the momentum, and this equation is the same as introduced in \eqref{eq:epdiff}. We note in particular that the momentum can be a distribution. The velocity field $v(\theta,t)$ lies in the space $\Vec(S^1)/\mathfrak{sl}_2(\R)$, where $\Vec(S^1)$ is the space of smooth vector fields on the circle and $\mathfrak{sl}_2(\R)$ is the Lie algebra of the group $\PSL_2(\R)$. However, the momentum $m(\th)\in\Hor$, where $\Hor=\{v=\sum \hat{v}_k e^{ik\th}:\hat{v}_k=0, k=0,\pm 1\}$. $\Hor$ is the orthogonal complement to $\mathfrak{sl}_2(\R)$ in $\Vec(S^1)$. We will refer to it as the {\em horizontal space}.

The energy is given by
\begin{align*}
  E(t) = \langle m(\cdot, t), v(\cdot, t) \rangle = \| v(\cdot, t) \|_{W P}^2,
\end{align*}
where $\langle \cdot, \cdot \rangle$ is the $L^2$ inner product on $[0, 2\pi]$. $E(t)$ is constant in time under geodesic flow. The $v \rightarrow m$ map may be inverted by the relation $v(\theta,t) = G*m(\theta,t)$, where $G$ is the Green's function $G(\th)$ of the $W P$ operator $L$. The Green's function $G(\th)$ is obtained as a solution to $L G = \mathrm{Proj}(\delta_0)$, where $\delta_0$ is the Dirac measure centered at $\theta = 0$, and $\mathrm{Proj}(\delta_0)$ is the projection of $\delta_0$ onto the horizontal space $\Hor$.
%, and $Proj$ is the projection of the Dirac function onto the horizontal space $\Hor$. 
 
The Green's function for the $W P$ operator \eqref{eq:wp-operator} is defined up to addition of any element of $\mathfrak{sl}_2(\R)$, $a+b\sin\th+c\cos\th$ for some constants $a,b,c$. We normalize $G(\th)$ as in \cite{kushnarev_teichons:_2009} so that it lies in the horizontal space $\Hor$:
\begin{align}\label{eq:greens-function}
  G(\theta) = \left(1 - \cos \theta\right) \log\left[ 2(1 - \cos \theta)\right] + \frac{3}{2} \cos \theta - 1.
\end{align}

\section{Teichons, Singular Solutions of EPDiff} \label{s:teichons}
The EPDiff equation \eqref{EPDiff} admits momenta solutions that, once initialized as a sum of $N$ Dirac measures, remain a sum of $N$ Dirac measures for all time \cite{camassa_integrable_1993, holm_momentum_2005}. In reference to this self-similarity property, these singular solutions are named teichons (or an $N$-teichon). In this paper, the velocity field that defines the geodesic between two shapes will be approximated by an $N$-teichon.

For a solution to EPDiff \eqref{EPDiff}, we employ the $N$-teichon ansatz
\begin{subequations}
\begin{align}\label{eq:m}
m(\theta,t) &= \sum_{j=1}^{N}p_j(t)\delta(\theta-q_j(t)),\\
\label{eq:v}
v(\theta,t) &= \sum_{j=1}^{N}p_j(t) G(\theta-q_j(t)),
\end{align}
\end{subequations}
where $\delta \triangleq \delta_0$ is the origin-centered Dirac mass. Plugging these expressions into EPDiff (\ref{EPDiff}), we obtain a system of ODEs describing the evolution of the momentum coefficients $p_k$ and the teichon locations $q_k$:
\begin{equation}
\left\{
\begin{array}{l}
\dot{p}_k = -p_k\sum_{j=1}^{N} p_j G'(q_k-q_j),\\*[.3cm]
\dot{q}_k = \sum_{j=1}^{N} p_j G(q_k-q_j).
\end{array}
\right.
\label{sol_ode}
\end{equation}
As mentioned in Section \ref{s:epdiff} the momenta $m$ lie in the horizontal space, i.e. $m(\theta,t)$ must have vanishing 0th and $\pm 1$st Fourier coefficients. Using \eqref{eq:m}, we obtain a set of three constraints for $(q_k,p_k)$, linear in $p_k$:
\begin{equation}
\sum_{j=1}^N p_j=\sum_{j=1}^N p_j e^{iq_j}=\sum_{j=1}^N p_j e^{-iq_j}=0. 
\label{wp_cond}
\end{equation}
If they are satisfied at time $t=0$ they will be satisfied for all $t$. The teichons never collide: i.e. the teichon locations $q_k$ retain their initial ordering on $S^1$ for all time. However, it is known that most initial configurations lead to an exponential decay of teichon separation, and an exponential increase in the momentum coefficients \cite{kushnarev_teichons:_2009}. This introduces numerical difficulties in solving \eqref{sol_ode}.

\section{Numerical method}
\label{s:nummethod}
Our numerical method has many components, and each of them requires some discussion. We proceed through these components as follows: we discuss construction of welding maps in Section \ref{sec:algorithm-welding}. The overarching shooting method is presented in Section \ref{sec:algorithm-shooting}, and the cross-ratio matching term is introduced in Section \ref{sec:algorithm-matching}. The construction of this matching term is nontrivial, involving a Delaunay triangulation of points in the complex plane, and Section \ref{sec:algorithm-delaunay} highlights these considerations. The (standard Euclidean) gradient of the matching term is necessary in order to update our initial teichon guess, and computation of the gradient is discussed in Section \ref{sec:algorithm-gradient}. It is well-known that the gradient is not the optimal search direction for optimization on Riemannian manifolds, and in addition the computed gradient does not satisfy the momentum admissibility conditions \eqref{wp_cond}; in Section \ref{sec:algorithm-optimization} we transform and project the gradient to address these concerns. Finally, Section \ref{sec:algorithm-initial-guess} discusses refinement strategies for obtaining good initial guesses, and Algorithm Listing \ref{alg:method} presents the full algorithm.

\subsection{Generating fingerprints}\label{sec:algorithm-welding}
We first consider the task of constructing fingerprints. In practice, motivated especially by computer vision, we will be given an ordered collection of points $\{z_m\}_{m=1}^M \subset \C$ lying on some closed curve in the complex plane. (This is our shape.) Naturally this does not define a continuous curve in the plane, but this situation is realistic in an application setting. 

From this discrete data, we use conformal welding to construct a continuous weld $\psi \in \PSL_2(\R) \backslash \Diff(S^1)$. This task requires only the ability to construct conformal maps between the unit disc and a region in the complex plane defined by the $z_m$. While many methods are suitable (in particular the methods described in \cite{sharon_2d-shape_2006}) we choose the Zipper algorithm \cite{marshall_convergence_2007}, which computes the map via a discretization of the Loewner differential equation. One main strength of the algorithm is its sequential nature -- computation of the entire map is a method that simply iterates over the index $m$ on the input points $z_m$ in an explicit way. In particular, the total work required is $\mathcal{O}(M^2)$. 

This should be contrasted with a popular competitor, numerical Schwarz-Christoffel mapping \cite{driscoll_schwarz-christoffel_2002,driscoll_algorithm_2005} requiring the solution of a size-$M$ optimization problem, which can exhibit convergence problems, which are exacerbated when the fingerprint has very large, or very small derivatives. The sequential nature of Zipper means that one computes the full conformal map as a composition of intermediate maps; when the fingerprint has exponentially large derivatives, such a compositional strategy is more computationally robust. Our experience indicates that Zipper algorithm is more efficient, and more resilient, for generating welding maps. We refer the reader to \cite{marshall_convergence_2007} for details on implementation of the Zipper algorithm, where convergence in the Hausdorff metric is proven. 

The situation of very large, or very small derivative values for a fingerprint is called ``crowding" in the literature: 
Recall that conformal welds given by \eqref{eq:weld-definition} are a composition of conformal maps $f_\inter$ and $f_\exter$. Given a conformal map $f$ whose domain is the unit disk $\D_\inter$, it is {\em crowded} if 
\begin{align*}
  R = \frac{\max_{x\in S^1} |f'(x)|}{\min_{x \in S^1} |f'(x)|}
\end{align*}
is large, where the definition of `large' depends on the finite-precision arithmetic being performed. A good rule-of-thumb for `large' is when $R$ is inversely proportional to machine precision, $R \sim \epsilon^{-1}_{\mathrm{mach}}$. A reasonable characterization of a crowded weld then is when the product of the $R$ values for the interior and exterior maps is on the order $\epsilon_{\mathrm{mach}}^{-1}$. Conformal welds $\psi$ are often stored as point-evaluations $(\theta_{\exter,m}, \theta_{\inter,m})$ that are computed from boundary values of the conformal maps $f_\inter$ and $f_\exter$. Any point-evaluations that are sampled in regions of $S^1$ where either map is crowded will coalesce to machine precision; that is,
\begin{align*}
  \frac{\theta_{\inter,m+1} - \theta_{\inter,m}}{|\theta_{\inter,m}|} \sim \epsilon_{\mathrm{mach}}, \hskip 20pt \textrm{or} \hskip 20pt
  \frac{\theta_{\exter,m+1} - \theta_{\exter,m}}{|\theta_{\exter,m}|} \sim \epsilon_{\mathrm{mach}}.
\end{align*}
When this happens, the weld $\psi$ is effectively not a diffeomorphism to machine precision, and it is difficult to accurately compute particle locations in crowded regimes. In the context of flow under EPDiff, particles $\theta_m$ may, at $t=0$, start in an uncrowded configuration (that is, the $t=0$ weld is not crowded), and then flow to a crowded configuration (that is, the $t=1$ weld is crowded). 
%In this case again a direct computational strategy will not be accurate: see subsection \ref{ss:other-shapes} for an illustration of this crowding.

Our algorithm does not ameliorate the underlying problem with crowding; indeed both Schwartz-Christoffel mapping and the Zipper algorithm do not produce accurate results for crowded shapes\footnote{We acknowledge the possibility that a solution is given in \cite{driscoll_numerical_1998}, but an application of this method to conformal welding is a separate, self-contained project in itself.}. However, we mention again that our experience is that Zipper is more robust when the shape is crowded.

Finally, we remark that although the input to the Zipper algorithm is a discrete set of points, the output is a continuous welding map $\psi$. Therefore in the sequel we continue to speak about continuous welding maps.

\subsection{Shooting method}\label{sec:algorithm-shooting}
We want to compute the geodesic between the two shapes, given by fingerprints $\psi_0$ and $\psi_1$. In other words we seek to find the velocity field $v$ corresponding to a geodesic for \eqref{EPDiff} such that the diffeomorphic evolution defined by \eqref{eq:v-evolution} satisfies the prescribed boundary conditions $\phi(t=0) = \psi_0$ and $\phi(t=1) = \psi_1$. We will solve this problem by approximating the velocity field by an evolving $N$-teichon solution to EPDiff.

Since \eqref{EPDiff} specifies the evolution of an initial velocity field, the goal then is to find the initial positions of teichons, $q_k(0)$, and initial teichon strengths, $p_k(0)$, such that the resulting velocity field will carry a template fingerprint $\psi_0$ as near as possible to the target $\psi_1$. The task of determining initial data to satisfy a two-point boundary value problem is well-studied and one of the more popular numerical methods to compute a solution is the shooting method \cite{press_numerical_2007}. Diffeomorphic matching in the context of shapes has also seen the recent application of shooting methods \cite{miller_geodesic_2006}. 

The idea of the shooting method is the following: start with an initial guess for the $q_k(0),p_k(0)$, construct the initial momentum $m(\th,0)=\sum p_k(0)G(\th-q_k(0))$, and then solve forward the equation (\ref{sol_ode}) to obtain time evolution of $q_k(t),p_k(t)$. This in turn will produce a time-varying velocity field $v(\th,t)$, which we integrate via the equation
\begin{subequations}
\begin{align}\label{eq:psi_ode}
v(\theta,t)&=\frac{\partial\phi}{\partial t}(\phi^{-1}(\theta,t),t),\\
\phi(\th,t=0)&=\psi_0.
\end{align}
\end{subequations}
to obtain $\phi(\th,t=1)$. The final computed fingerprint $\phi(\th,t=1)$, is compared with the target fingerprint, $\psi_1$. Based on this comparison, we modify the initial shot configuration $\{p_k(0), q_k(0)\}_{k=1}^N$ and repeat the process. Because the $p_k$ must satisfy constraints that depend on $q_k$, our shooting method only changes the teichon momenta $p_k$; changing $q_k$ is certainly possible but requires admissibility constraints whose application is more involved. We have found that varying only the momenta allows us to represent a large variety of shapes.

\begin{figure}
  \begin{center}
    \resizebox{1.0\textwidth}{!}{
      \includegraphics{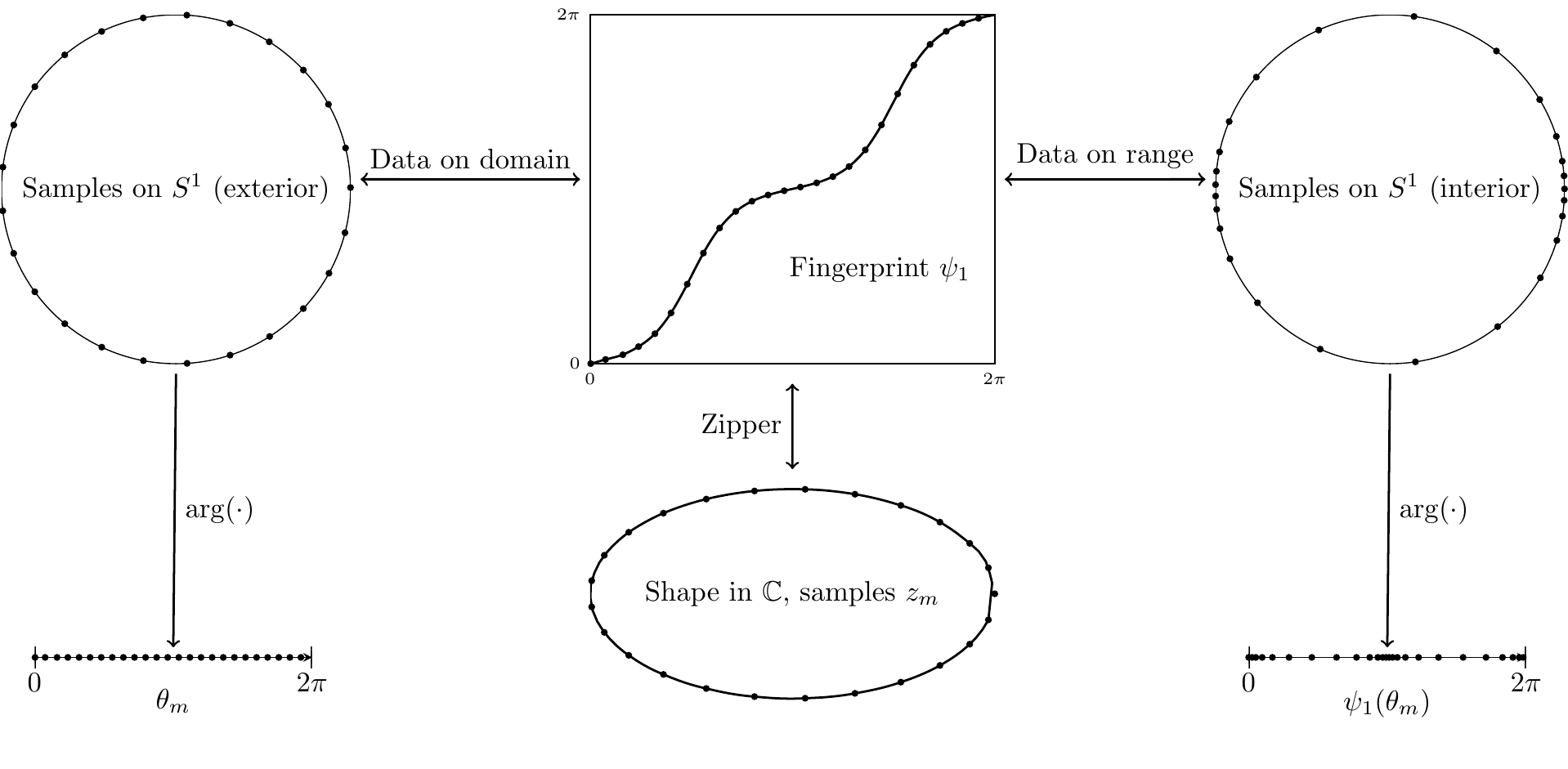}
    }
  \end{center}
  \caption{For a terminal shape, relationship between the desired landmark positions $\psi_1(\theta_m)$, and the welding map generated from the shape samples $z_m$. The initial landmark locations $\alpha_m = \psi_0(\theta_m)$ are not shown.}
  \label{fig:welding-transformation}
\end{figure}

In practice, we cannot match $\psi_1$ up to infinite precision, so we resort to inexact matching via some discrete set of landmark points. Given the discussion from Section \ref{sec:algorithm-welding}, it is sensible to choose $M$ landmarks corresponding to the images of the terminal shape samples $z_m$. (Here, {\em terminal} means the target shape at $t=1$.) We recall that fingerprints $\psi$ are defined through conformal maps $f$: $\psi_0 = (f_{0,\inter}^{-1} \circ f_{0,\exter})|_{S^1}$, and similarly for $\psi_1$. We consider landmark points $\theta_m \in [0, 2\pi)$ on the exterior defined by 
\begin{align}\label{eq:thetam-definition}
  \exp (i \theta_m) \triangleq f^{-1}_{1,\exter}(z_m).
\end{align} 
We track these landmarks as they flow from the initial shape: let $\alpha_m = \psi_0(\theta_m) =  \arg \circ f^{-1}_{0,\inter} \circ f_{0,\exter}(z_m)$. We flow these landmarks using \eqref{eq:psi_ode} to $t=1$ and in principle we wish to compare their locations with the exact terminal locations $\psi_1(\theta_m) = f^{-1}_{1,\inter}(z_m)$. The $t=1$ images of $\alpha_m$ under an $N$-teichon evolution are determined by 
\begin{equation}\label{eq:alpha-evolution}
  \begin{aligned}
    \dot{\alpha}_m &= \sum_{n=1}^N p_n G(\alpha_m - q_n), \\
    \alpha_m(0) &= \psi_0(\theta_m).
  \end{aligned}
\end{equation}
The particular choice that $\alpha_m(0) = \psi_0(\theta_m)$ is not the only choice one could make, and we do not claim it is optimal; however, our results indicate that such a choice performs quite well in many situations. The choice of landmark locations $\alpha_m$ and the $t=0$ teichon configuration $q_k$ need not be related. With these $M$ landmark locations $\alpha_m(1)$, we must compare fingerprints at $t=1$. We thus need a matching functional $E_2(\psi_1,\psi(\theta,1))$ that reflects the closeness of the fingerprints. The relation between the terminal shape samples $z_m$ and the landmarks $\theta_m$ is shown in Figure \ref{fig:welding-transformation}.

%\tikzsetnextfilename{welding-transformation}

\subsection{Matching fingerprints}\label{sec:algorithm-matching}

The shooting method relies on computation of a matching term, which we call $E$. This matching function compares the fingerprint computed with an $N$-teichon evolution with a target fingerprint. Note that the standard type of matching functional for geodesic shooting in this context would take the form 
\begin{align}\label{eq:full-matching-term}
  E &= \int_0^1 \| v(\cdot, t) \|_{W P}^2 \dx{t} + \lambda \sum_{m=1}^M \left(\psi_1(\theta_m) - \alpha_m(1) \right)^2\\
\notag
&= E_1 + E_2.
\end{align}
for the landmark choices $\{\theta_m\}_{m=1}^M \subset S^1$, and a scalar $\lambda \geq 0$. The weight $\lambda$ defines the relative importance between the first energy term, and the second matching term. A shooting procedure would iterate on the initial data in an attempt to minimize this functional.
Our version employs two variations. First, we are not minimizing energy: geodesics on $\uts$ are unique \cite{gay-balmaz_geometry_2009} so that the value of the energy is irrelevant since only one path exists between $\psi_0$ and $\psi_1$. Therefore we entirely omit the first energy term $E_1$ in \eqref{eq:full-matching-term}.

\begin{figure}
  \begin{center}
    \resizebox{1.0\textwidth}{!}{
      \includegraphics{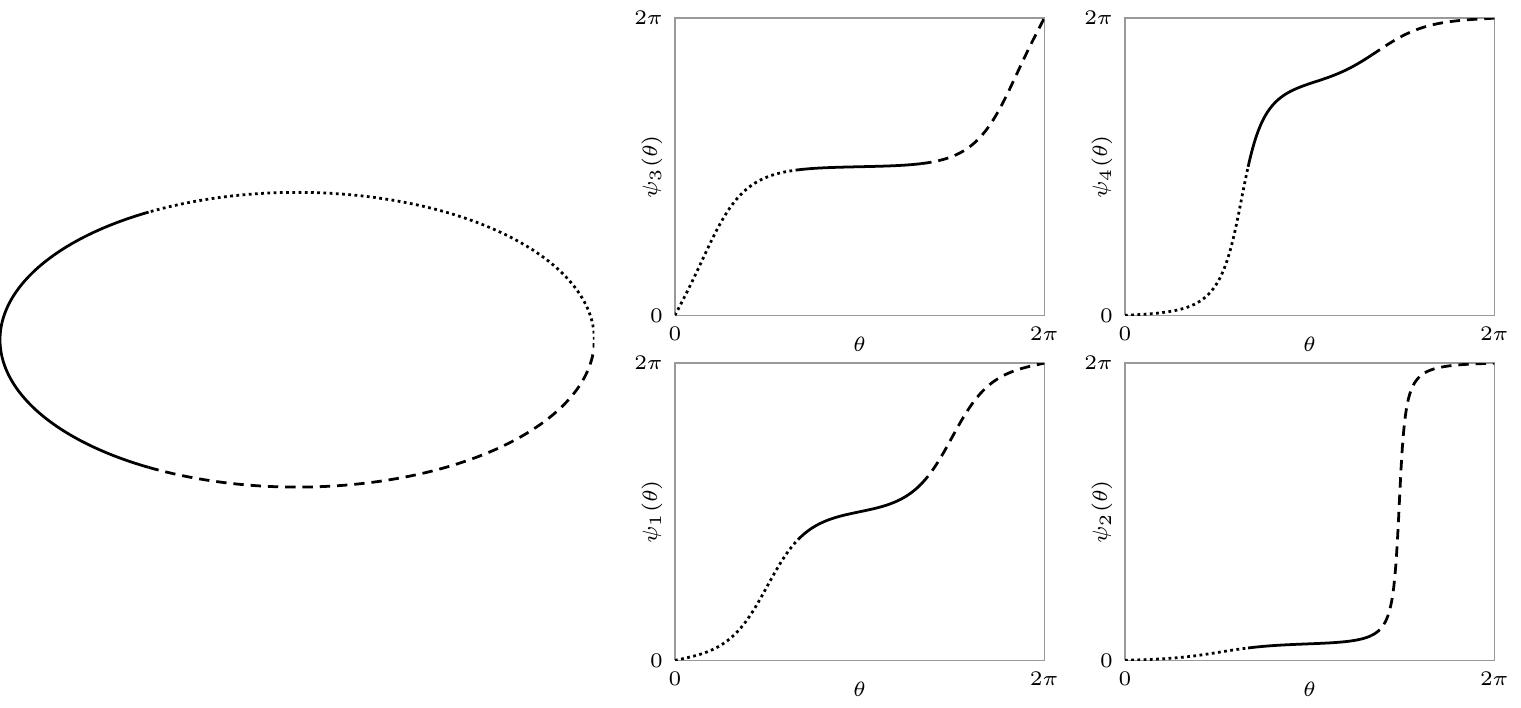}
    }
  \end{center}
  \caption{Left: an ellipse. Right: various welding map representatives of the ellipse shape corresponding to the same equivalence class in $\uts$.}
  \label{fig:welding-maps}
\end{figure}

Our method also employs a different landmark matching term. We cannot directly employ an $\ell^2$-type distance as given in \eqref{eq:full-matching-term} because of the M\"obius invariance of welding maps. For comparison, we show four different welding maps associated with the same shape in Figure \ref{fig:welding-maps}. We want any matching term that we devise to assign zero distance between any pair of welding maps in the figure. However, using an $\ell^2$ type distance to compare them is clearly misleading. In particular, given any fingerprint $\psi_1$ and $\varepsilon > 0$, there is a M\"obius map $A_m \in PSL_2(\R)$ such that the pointwise distance between $\psi_1$ and $A_m \circ \psi_1$ at $\theta_m$ is within $\varepsilon$ of the maximum matching distance:
\begin{align*}
  (\psi_1(\theta_m) - \left(A_m \circ \psi_1\right)(\theta_m))^2 \geq 4 \pi^2 - \varepsilon.
\end{align*}
%\tikzsetnextfilename{welding-maps}
Given $\psi_0$ and $\psi_1$, we cannot determine the appropriate M\"obius maps so that $\psi_0$ and $\psi_1$ have the ``same" normalization until we have already computed the geodesic connecting them. 
Since the matching term should be insensitive to self-maps of the disk, we amend the landmark matching term to be robust with respect to projective transformations on $\C$. The most general projective invariant quantity of a 4-tuple of points in the complex plane is the cross-ratio and forms the basis for our matching term.

Let $(z_1,z_2,z_3,z_4) \subset \C$ be a 4-tuple of points. The cross-ratio of these four points is defined by 
\begin{align}\label{eq:cross-ratio}
C(z_1,z_2,z_3,z_4) = \frac{(z_1-z_3)(z_2-z_4)}{(z_2-z_3)(z_1-z_4)}.
\end{align}
$C$ is invariant under the M\"obius transformations: $C(z_1,z_2,z_3,z_4)=C(A(z_1),A(z_2),A(z_3),A(z_4))$ for any M\"obius map $A$. One can also show that if $z_j \in S^1$ for all $j$, then $C \in \R$. 

The invariance of the cross-ratio under M\"obius transformations will allow us to compare fingerprints that have different M\"obius normalizations. After evolving according to EPDiff, we have $M$ landmark locations on $S^1$ that specify pre-images of shape vertices under $\phi_{\mathrm{int}}$. This suggests that we can only resolve the shape up to these $M$ vertices, and that furthermore we can only use cross-ratios of these pre-images in order to have the M\"obius invariance. 

A method to uniquely encode information about a polygon with $M$ vertices has been proposed in \cite{driscoll_numerical_1998}. The basic idea is that a polygon is uniquely identifiable if the vertex angles of a polygon are specified in conjunction with $M-3$ carefully chosen cross-ratios of quadrilaterals. These quadrilaterals are constructed from the {\em Delaunay triangulation} of the polygon.

\subsection{Delaunay triangulation}\label{sec:algorithm-delaunay}
Let $P$ be a simple polygon. A {\it triangulation} of $P$ is a division of $P$ into non-degenerate triangles whose vertices are vertices of $P$. The triangles intersect only at a vertex or at an entire edge. A {\it Delaunay triangulation} of polygon $P$ is a triangulation such that no point in $P$ is inside the circumcircle of any triangle in the Delaunay triangulation. Triangle edges of the triangulation that are not polygon edges are called {\it diagonals}.

It is known \cite{marshall_bern_mesh_1992} that every $P$ has at least one Delaunay triangulation with the following property. If $d$ is the diagonal, let $Q(d)$ be the quadrilateral, composed of the union of two triangles on either side of $d$. Then the sum of two opposite interior angles of $Q(d)$ that are split by $d$ is at least $\pi$. A Delaunay triangulation of an $n$-polygon $P$ can be computed in $O(n^2)$ steps. In our implementation we have used the MATLAB function \texttt{delaunay}.
%As has been shown in \cite{DrVa} if we choose cross-ratios through a Delauney triangulation, we are guaranteed
%Denote by $Z_j$ a four-tuple of points on a circle. The matching functional is the sum of squared differences of the cross-ratios:

\begin{figure}[h]
\begin{center}
\includegraphics[width=5cm]{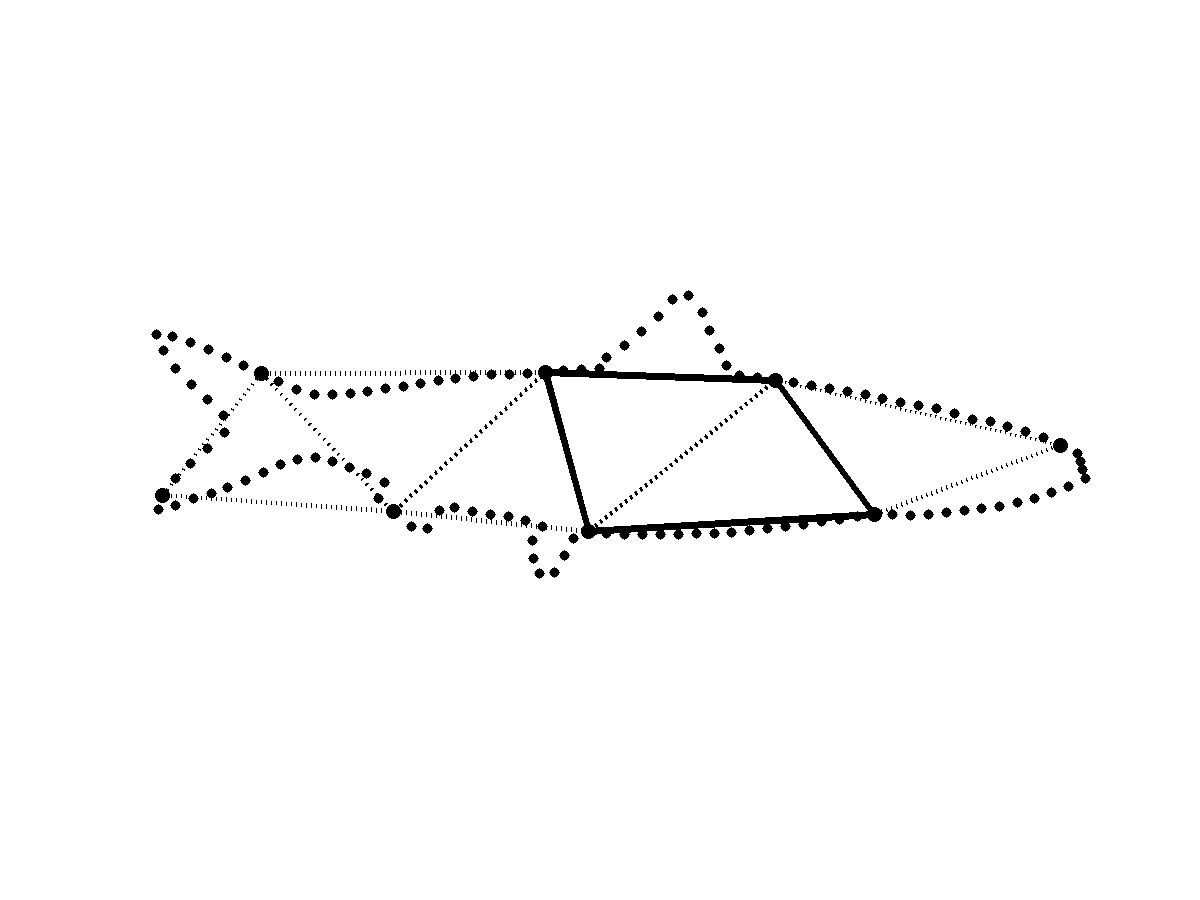}
\includegraphics[width=5cm]{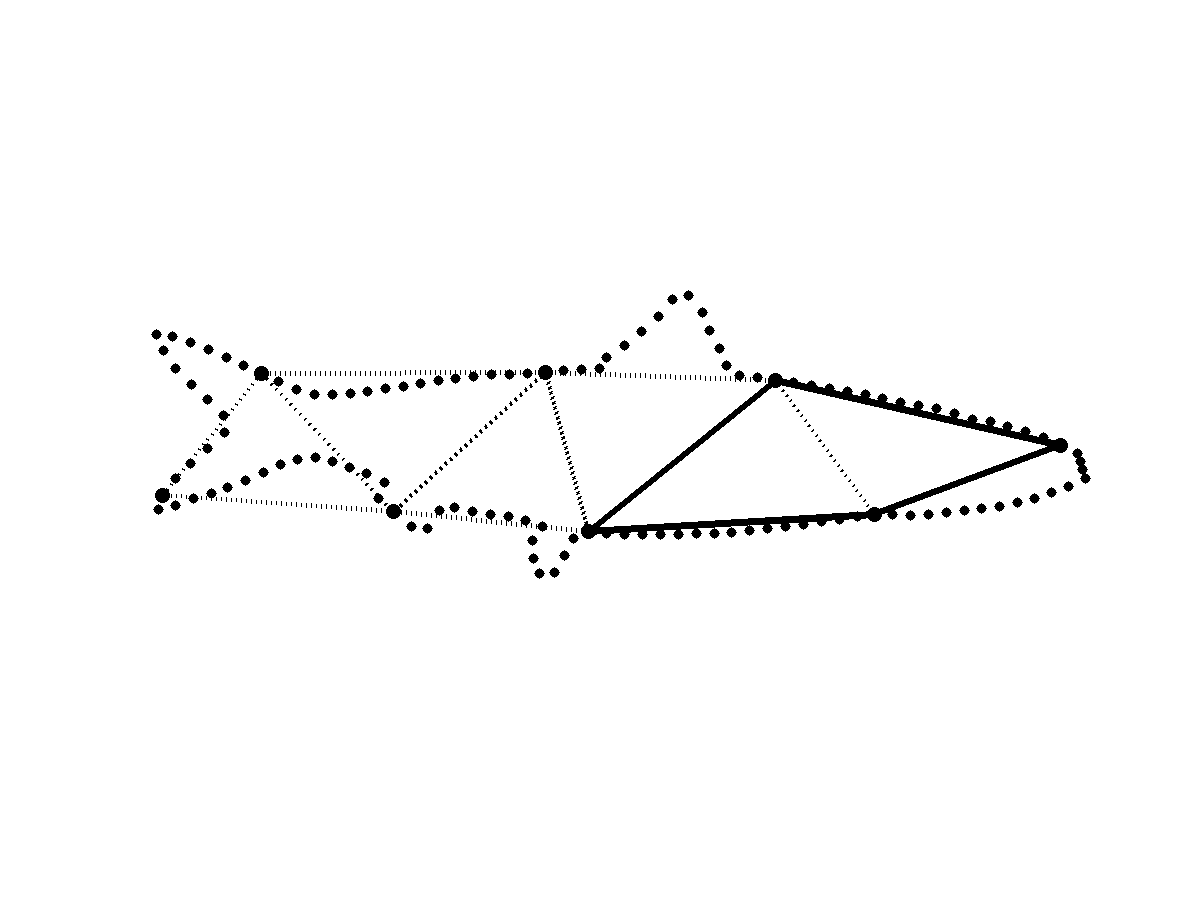}
\includegraphics[width=5cm]{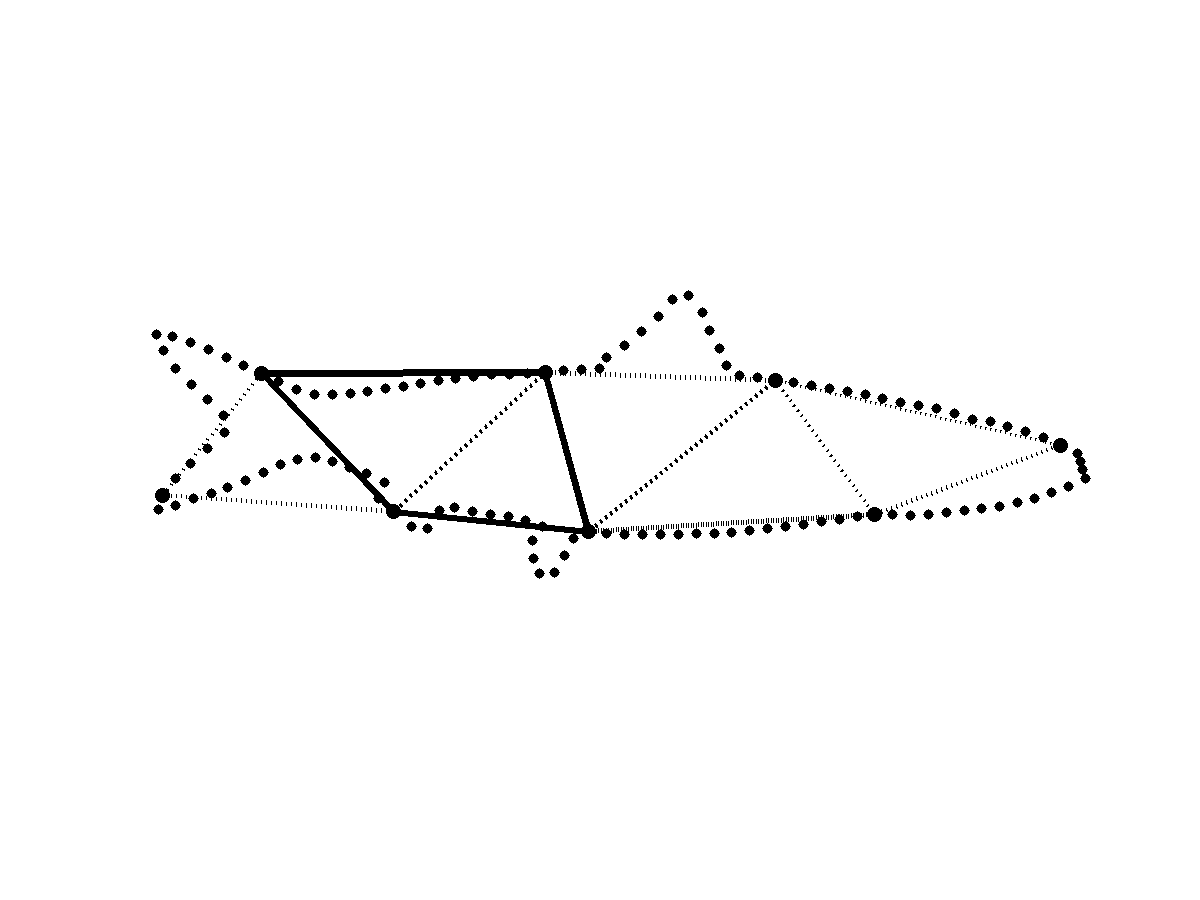}
\includegraphics[width=5cm]{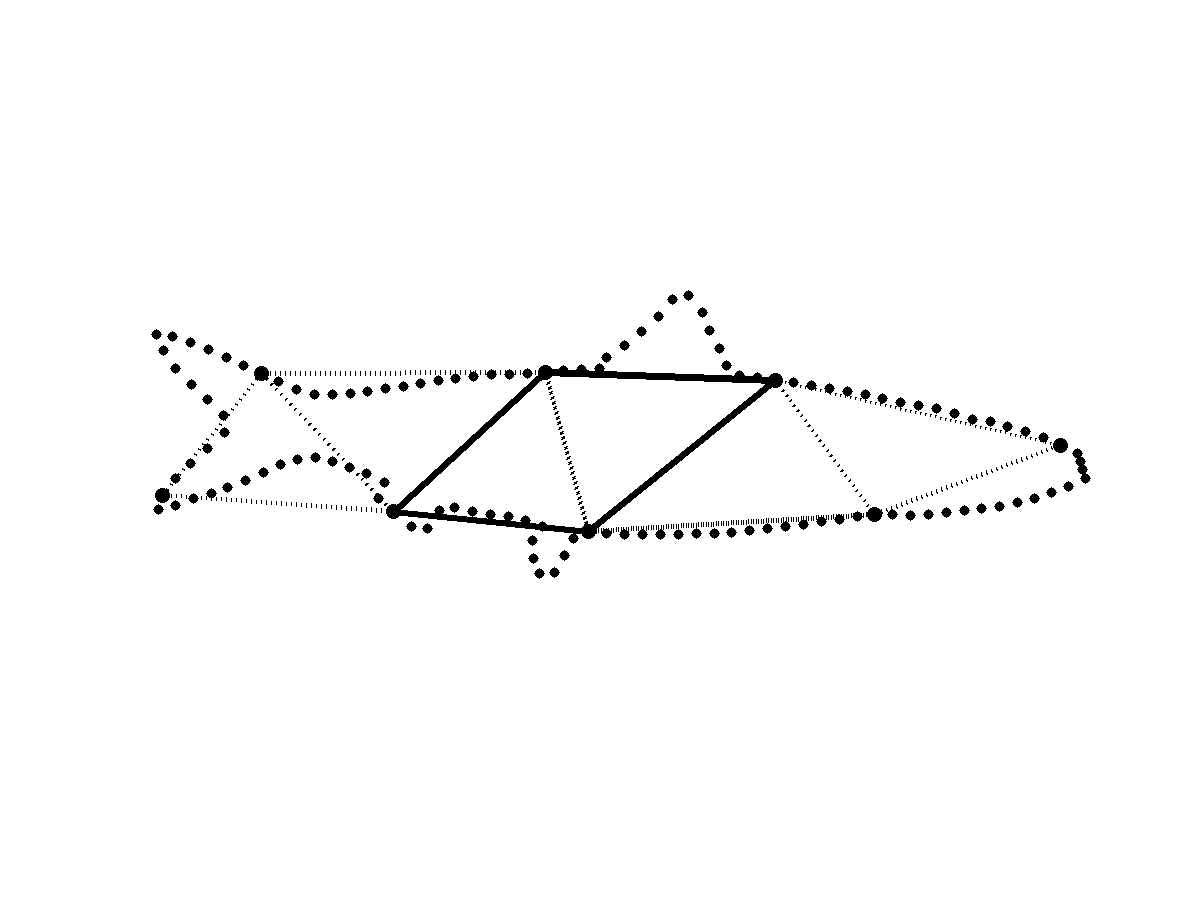}
\includegraphics[width=5cm]{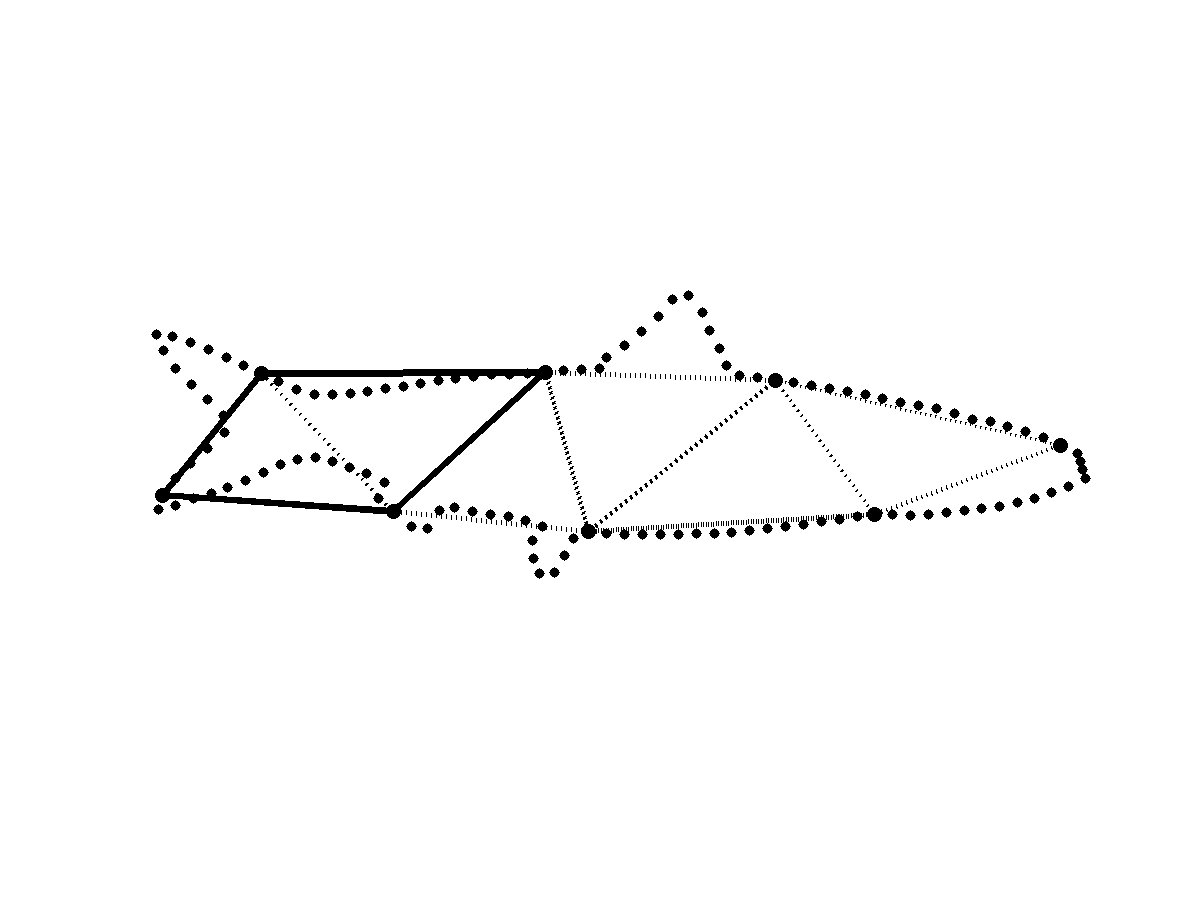}
\end{center}
\caption{The 4-tuples used in computing of cross ratios in the matching functional (\ref{eq:cr}) are connected by solid lines. Dotted lines depict Delaunay triangulation on the coarsest subset of $N$ points representing the fish ($N=8$). For $N$ points used in computing the cross ratios we get $N-3$ four-tuples.}
\label{fig:dt-quads}
\end{figure}

It is well known that any $n$-vertex simple polygon $P$ has a triangulation consisting of $n-2$ triangles. In addition it has exactly $n-3$ distinct diagonals. The vertices of a quadrilateral $Q(d)$ associated with each diagonal are used in the computation of the cross-ratios, see Figure \ref{fig:dt-quads}. More specifically, it is shown in \cite{driscoll_numerical_1998} that using cross-ratios computed using this choice of quadrilaterals uniquely characterizes the original polygon $P$.

With this in mind, we let $\{i_{k,1}, i_{k,2}, i_{k,3}, i_{k,4} \}_{k=1}^K$ denote the $K$ 4-tuples of cross-ratios from \cite{driscoll_numerical_1998}. If we have $M$ points for our shape, $i_{k,j}$ is an index taking a value in $1, \ldots, M$ indicating which point to use in the $k$th cross-ratio. We take the matching functional to be the relative error of the discrete $\ell^2$ cross-ratio difference:
\begin{align}
\label{eq:cr}
%E(\phi_1,\phi(\theta,1)) = \sum_{k=1}^K \left(\frac{C^{\phi_1}(Z_j)-C^{\phi(\theta,1)}(Z_j)}{C^{\phi_1}(Z_j)}\right)^2.
E_2(\psi_1(\theta),\phi(\theta,1)) = \frac{1}{K} \sum_{k=1}^K \left(1 -
                  \frac{[C\circ e \circ \phi(\cdot,1)] (\theta_{i_{k,1}}, \theta_{i_{k,2}}, \theta_{i_{k,3}}, \theta_{i_{k,4}})}
                        {[C\circ e \circ \psi_1] (\theta_{i_{k,1}},\theta_{i_{k,2}}, \theta_{i_{k,3}}, \theta_{i_{k,4}})} \right)^2,
\end{align}
where $e(x) \triangleq \exp(i x)$ is the complex exponential, and the $\theta_m$ are the points defined by \eqref{eq:thetam-definition}.

We emphasize that the construction \eqref{eq:cr} for the matching term is automated: the identification of the required quadrilaterals making up the cross-ratio term is automatically computed using the Delaunay triangulation of the shape. After a one-time run of the Delaunay triangulation, the formula \eqref{eq:cr} is also explicit: the indices $i_{k,j}$ are known and stored.

\subsection{Gradient with respect to $p(0)$}\label{sec:algorithm-gradient}
In order to adjust the initial momenta to reach the target, we need to compute the gradient of the matching functional, $E$, with respect to $p(0)$. Let $p_0$ denote $p(0)$, then the gradient of the energy is given by
\begin{align}\label{eq:E-chain-rule}
  \frac{\partial E}{\partial p_0} = \frac{\partial E}{\partial\alpha}\frac{\partial \alpha}{\partial p_0},
\end{align}
where $\alpha$ is the length-$M$ vector of landmarks solving \eqref{eq:alpha-evolution} at time $t=1$. Define $\beta_m \triangleq \pfpx{\alpha_m}{p_0}$, $\pi_n \triangleq \pfpx{p_n}{p_0}$, and $\chi_n \triangleq \pfpx{q_n}{p_0}$, each of which is a $1 \times N$ vector for $m = 1, \ldots, M$ and $n = 1, \ldots, N$. These parameters may be computed by a system derived from \eqref{sol_ode} and \eqref{eq:alpha-evolution}:
\begin{subequations}
\label{eq:modified-ode}
\begin{align}
  \dfdx{\beta_m}{t} &= \beta_m \sum_{k=1}^N \left[ \pi_k G(\alpha_m - q_k) + p_k G'(\alpha_m - q_k) (\beta_m - \chi_k) \right] \\
  \dfdx{\pi_n}{t} &= -\pi_n \sum_{k=1}^N p_k G'(q_n - q_k) - p_n \sum_{k=1}^N \left[ \pi_k G'(q_n - q_k)  + p_k G''(q_n - q_k) (\chi_n - \chi_k)\right] \\
  \dfdx{\chi_n}{t} &= \sum_{k=1}^N \left[ \pi_k G(q_n - q_k) + p_n G'(q_n - q_k) (\chi_n - \chi_k) \right],
\end{align}
\end{subequations}
where we assign $G''(0) = 0$. The full system \eqref{eq:modified-ode}, \eqref{eq:alpha-evolution}, and \eqref{sol_ode} can be solved in parallel to determine $\beta_m(1) = \pfpx{\alpha}{p_0}$ to be used in \eqref{eq:E-chain-rule}. One can explicitly compute $\pfpx{E}{\alpha}$ from \eqref{eq:cross-ratio}, and the definition of the complex exponential $e(\cdot)$. Thus the energy gradient \eqref{eq:E-chain-rule} is computable.

\subsection{Optimization with the gradient}\label{sec:algorithm-optimization}

There are three tasks yet to be accomplished before we can update the initial momentum distribution: 
\begin{itemize}
  \item $\partial E/\partial p_0$ is not an admissible momentum distribution, so we must project it into the appropriate space
  \item on non-Euclidean Riemannian manifolds the gradient does not point in the direction of steepest ascent; we require the {\em natural gradient}
  \item gradient descent is the most basic of optimization methods; we employ a nonlinear conjugate gradient update to accelerate convergence
\end{itemize}
This subsection discusses these considerations.

\subsubsection{Projecting the gradient}\label{sec:algorithm-gradient-projection}
For the Teichon evolution system (\ref{sol_ode}) to be valid one needs to have a bijection between the Lie algebra $\g$ and its dual $\g^*$. This bijection is provided by the Weil-Petersson operator
$L$ and it's inverse, convolution with the Green's function. These operators are a bijection only on the horizontal space $\Hor=\{w\in\Vec(S^1):\hat{w}_k=0, k=0,\pm 1 \textrm{ and } \| w \|_{W P} < \infty\}$, where the $\hat{w}_k$ are Fourier coefficients of the periodic function $w(\th)$ defined on the circle. 
In other words:
\begin{align*}
L&:Hor\rightarrow Hor^*,\\
G&:Hor^*\rightarrow Hor^,\\
L&G = G*L\delta = Proj(\delta).
\end{align*}
Here, $\mathrm{Proj}(\delta)$ is the projection of the delta function onto the horizontal subspace $Hor$. Therefore any updates we perform to the initial momentum distribution must happen on this space. In order for a momentum field $m \in \Vec(S^1)$ to lie in $Hor^\ast$, it must likewise have vanishing $0, \pm 1$ Fourier coefficients. This defines the three constraints given by \eqref{wp_cond} as discussed in Section \ref{s:teichons}.

While our starting guess for $p(0)$ will satisfy the conditions \eqref{wp_cond}, there is no guarantee that the update vector $\pfpx{E}{p_0}$ from \eqref{eq:p0-update} will satisfy those constraints. We must therefore obtain an element from $Hor^\ast$ given $\pfpx{E}{p_0}$.

Let $\Delta p^{new} = \pfpx{E}{p_0}$, which does not represent a member of $Hor^\ast$. We project this update vector into the space admissible updates: those that satisfy \eqref{wp_cond}. We proceed by computing the $WP$-closest member of $Hor^\ast$ to $\Delta p^{new}$. Let the $3 \times N$ matrix $F$ define the admissibility constraints:
\begin{align*}
F(q(0)) = 
\left(
\begin{matrix}
1 & \ldots &1\\
\cos q_1(0) &\ldots &\cos q_N(0)\\
\sin q_1(0) &\ldots &\sin q_N(0)
\end{matrix}
\right)
\end{align*}
%An $N$-vector $\Delta p$ representing a member $Hor^\ast$ lies in the kernel of $F(q(0))$. 
We wish to find an update vector $\widetilde{\Delta p}$ satisfying 
\begin{align}\label{eq:horizontal-minimization}
  \widetilde{\Delta p} = \argmin_{F(q(0)) \Delta p = 0} \| \Delta p - \Delta p^{new} \|_{W P^\ast},
\end{align}
where $\|\cdot\|_{WP^\ast}$ is the norm on $\Hor^\ast$ induced by the norm $\|\cdot\|_{W P}$ on $\Hor$:
\begin{align*}
  \|\Delta p\|_{WP^\ast} \triangleq \| G \ast \Delta p \|_{W P}.
\end{align*}
Let $G(q(0))$ be the $N \times N$ Gram matrix for $\{\mathrm{Proj}(\delta_{q_n(0)})\}_{n=1}^N$, whose entries are given by $G_{i, j}(q(0)) = G(q_i(0) - q_j(0))$, where $G(\cdot)$ is the Green's function \eqref{eq:greens-function}. Then \eqref{eq:horizontal-minimization} can be written as minimization of a quadratic objective subject to a linear constraint. The solution is a vector $\widetilde{\Delta p}$ given by 
\begin{align}\label{eq:gradient-projection}
\widetilde{\Delta p} = \Bigl[I - G^{-1}F^T(FG^{-1}F^T)^{-1}F\Bigr]\Delta p,
\end{align}
where above $G := G(q(0))$ and $F := F(q(0))$.

\subsubsection{The natural gradient}\label{sec:algorithm-natural-gradient}
It is well-known that performing gradient descent on Riemannian manifolds with the standard gradient is not the optimal gradient update. It is much more effective to use the natural gradient (see, e.g., \cite{amari_natural_1998}). The idea behind the natural gradient is the following: given an $N$-teichon $p$ on the space $Hor^\ast$, suppose $\widetilde{\Delta p}$ is the standard gradient direction (satisfying the admissibility constraints from the previous section). We wish to move in the direction that decreases the objective most, and so we must solve the problem
\begin{align*}
  \widehat{\Delta p} = \argmax_{\|p \|_{WP^\ast} = 1} \langle p, \widetilde{\Delta p} \rangle_{W P^\ast} = \argmax_{\|p \|_{WP^\ast} = 1} p^T G \widetilde{\Delta p}
\end{align*}
where the constraint $\| p \|_{WP^\ast} = 1$ may be written as $p^T G p = 1$ where $G$ is the Gram matrix for the $N$-teichon configuration of $p$. The solution to this problem is given by 
\begin{align}\label{eq:natural-gradient}
 \widehat{\Delta p} = \lambda G \Delta \tilde{p},
\end{align}
for some normalizing constant $\lambda$, which we hereafter set to unity. The new gradient direction $\widehat{\Delta p}$ is called the natural gradient. Combining this with \eqref{eq:gradient-projection}, the full update vector $\Delta \hat{p}$ given the unconstrained gradient is
\begin{align}\label{eq:proper-gradient}
  \widehat{\Delta p} = \Bigl[G - F^T(FG^{-1}F^T)^{-1}F\Bigr] \pfpx{E}{p_0}
\end{align}
In the sequel we will refer to the momentum update \eqref{eq:proper-gradient} that is both admissible (satisfying \eqref{wp_cond}) and natural (given by \eqref{eq:natural-gradient}) as the {\em proper} gradient.
%Then the new update to the initial momenta at each iteration is given by 
%\begin{align*}
%  p(0) \gets p(0) - \varepsilon \Delta \hat{p},
%\end{align*}
%for some $\varepsilon > 0$.

\subsubsection{Updating the shooting direction}\label{sec:algorithm-cg}
With the proper gradient vector $\widehat{\Delta p}$ given by \eqref{eq:proper-gradient}, we can proceed with standard optimization methods. The most straightforward is steepest descent: The update for the vector $p(0)$ at each iteration is given by
\begin{align}\label{eq:p0-update}
  p(0) \gets p(0) - \varepsilon \widehat{\Delta p}
\end{align}
where the choice $\varepsilon$ determines how far along the gradient direction $\widehat{\Delta p}$ we update. 

Convergence with gradient descent often stagnates when the iterative path taken by solutions follows a narrow valley; in such cases more sophisticated methods are required to render the iteration computationally efficient. A nonlinear conjugate gradient method is one such alternative. We employ the Polak-Ribi\`ere method, which makes use of the update vectors from the previous iteration. Let $\rho_n$ denote the proper gradient \eqref{eq:proper-gradient} at iteration $n$. Compute
\begin{align}\label{eq:cg-beta}
  \beta = \max\left\{0, \frac{\langle \rho_n, \rho_n - \rho_{n-1}\rangle_{W P^\ast}}{\langle \rho_{n-1}, \rho_{n-1} \rangle_{W P^\ast}} \right\}
\end{align}
Set the update direction at iteration $n$ to $\omega_n = \rho_n + \beta \omega_{n-1}$. Perform the update $p(0) \gets p(0) - \varepsilon^\ast \omega_n$, where $\varepsilon^\ast$ satisfies
\begin{align}\label{eq:cg-update}
  \varepsilon = \argmin_{\varepsilon > 0} E(p(0) - \varepsilon \omega_n)
\end{align}
To begin the iteration process, the first iteration is performed as a standard gradient descent update.

We utilize the standard nonlinear optimization tricks: standard gradient descent is performed at the initial stages to iterate close to a basin of attraction; then nonlinear conjugate gradient is employed to quickly converge.

\subsection{Initial configuration guesses}\label{sec:algorithm-initial-guess}
In all of our experiments, we specify the number of teichons $N$, which is fixed throughout iteration. We also determine the initial teichon configuration $q(0)$ by spacing the $N$ teichons equidistantly on $S^1$. This is not the optimal choice, especially for shapes with important features concentrated in a particular location, but our tests indicate that this is not a bad choice for many non-crowded shapes.

In taking the initial choice for $p(0)$ to be the zero vector, we have found that convergence takes an inordinate amount of time, or the iteration stagnates and convergence is not observed at all. One well-known way to combat this is to choose a better initial guess. We do this by first solving the minimization problem on a subset of the full collection of cross-ratios. We identify this subset by choosing cross-ratios that resolve coarse features of the shape.

We therefore implement a coarse-to-fine approach in the matching functional \eqref{eq:cr}. We compute the Delaunay triangulation on a coarse subset of the vertices of a given polygon $P$. For example, if a polygon $P$ has $M=128$ vertices, we compute a triangulation on dyadic subsets consisting of 8, 16, 32 and 64 points, and finally the full set of 128 points. We use a zero initial guess for $p(0)$ and run the minimization algorithm for the cross-ratios identified by the 8-point Delaunay Triangulation. The solution of the 8-point problem is used as the initial guess for the 16-point problem, and so on. Assuming a nested choice of points (i.e.~the 8 points are a proper subset of the 16 points, etc.) then we progressively build up the choice of cross-ratios until we utilize those for the full set of 128 points.

\begin{figure}[h]
\begin{center}
\includegraphics[width=5cm]{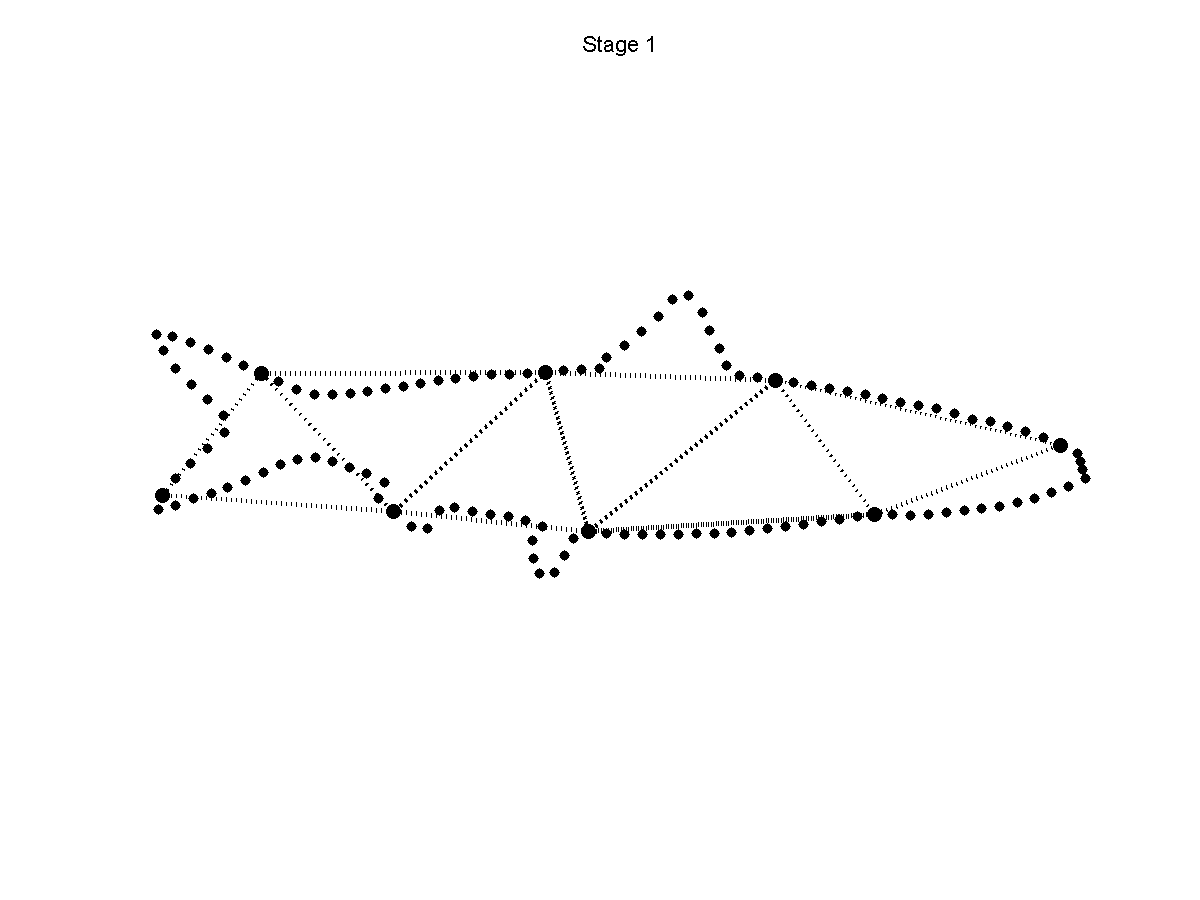}
\includegraphics[width=5cm]{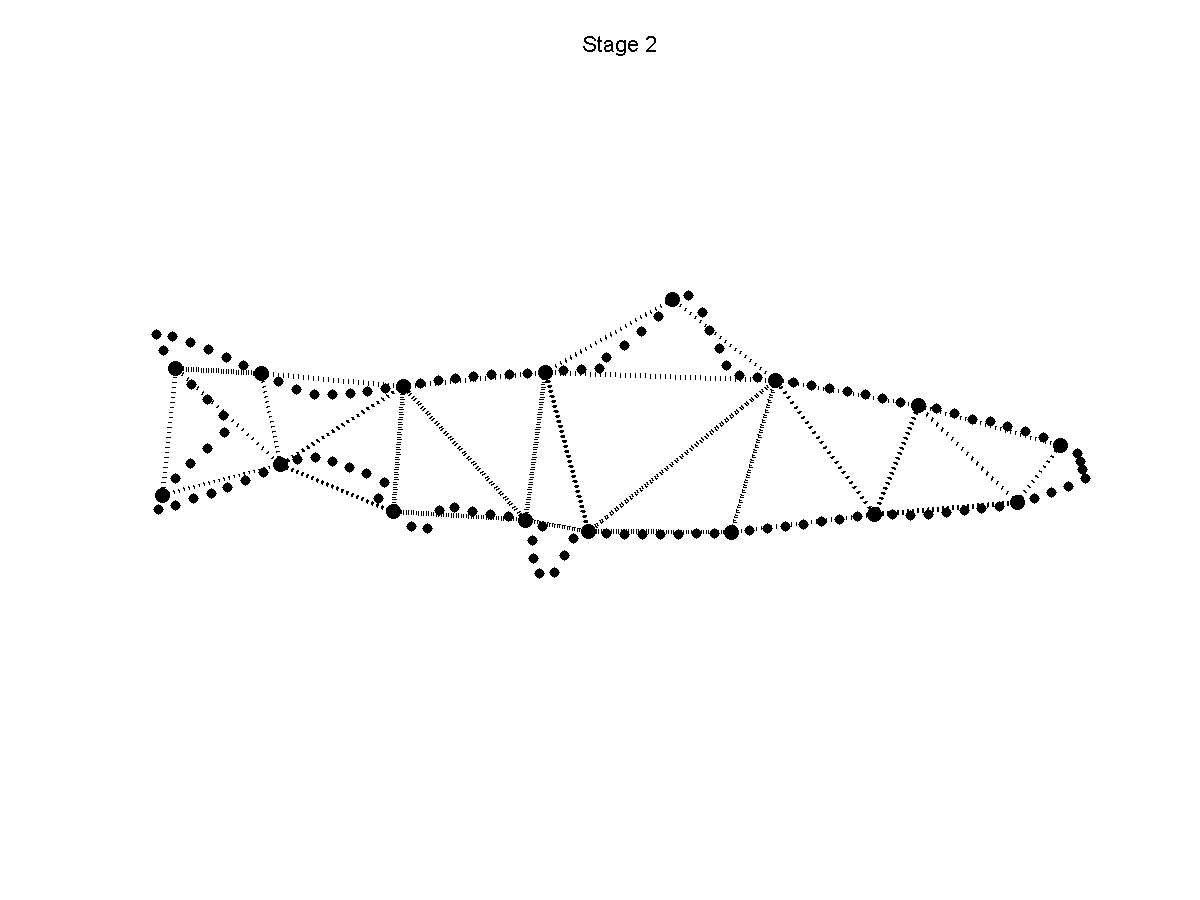}
\includegraphics[width=5cm]{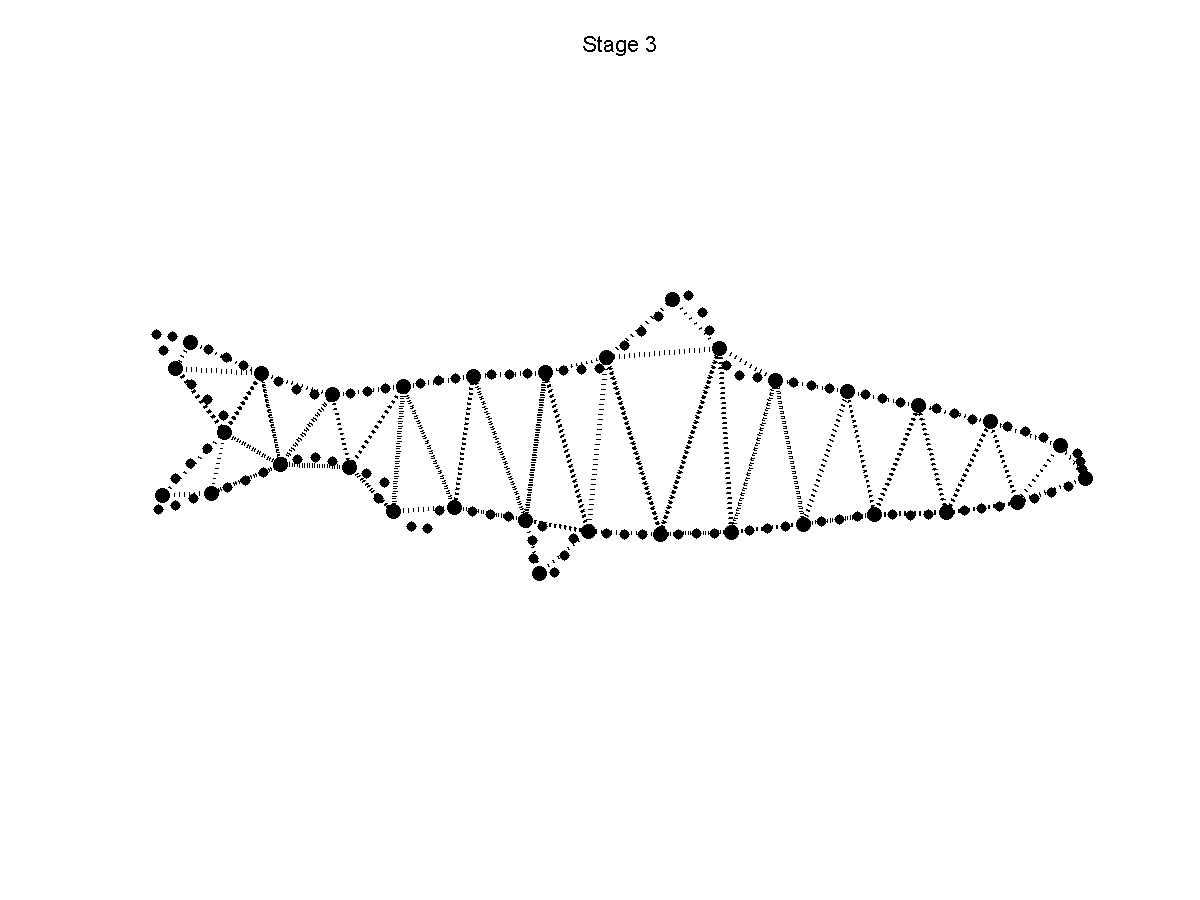}
\includegraphics[width=5cm]{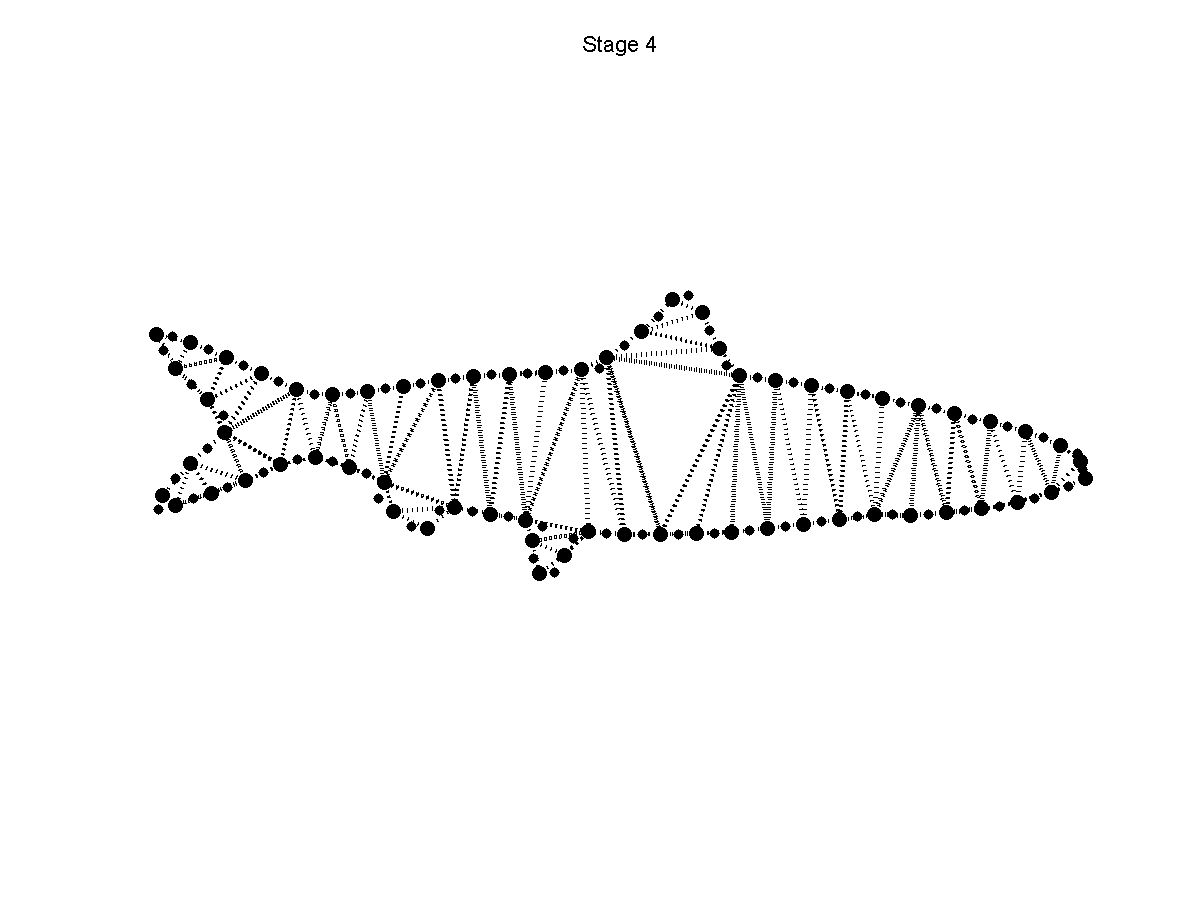}
\includegraphics[width=5cm]{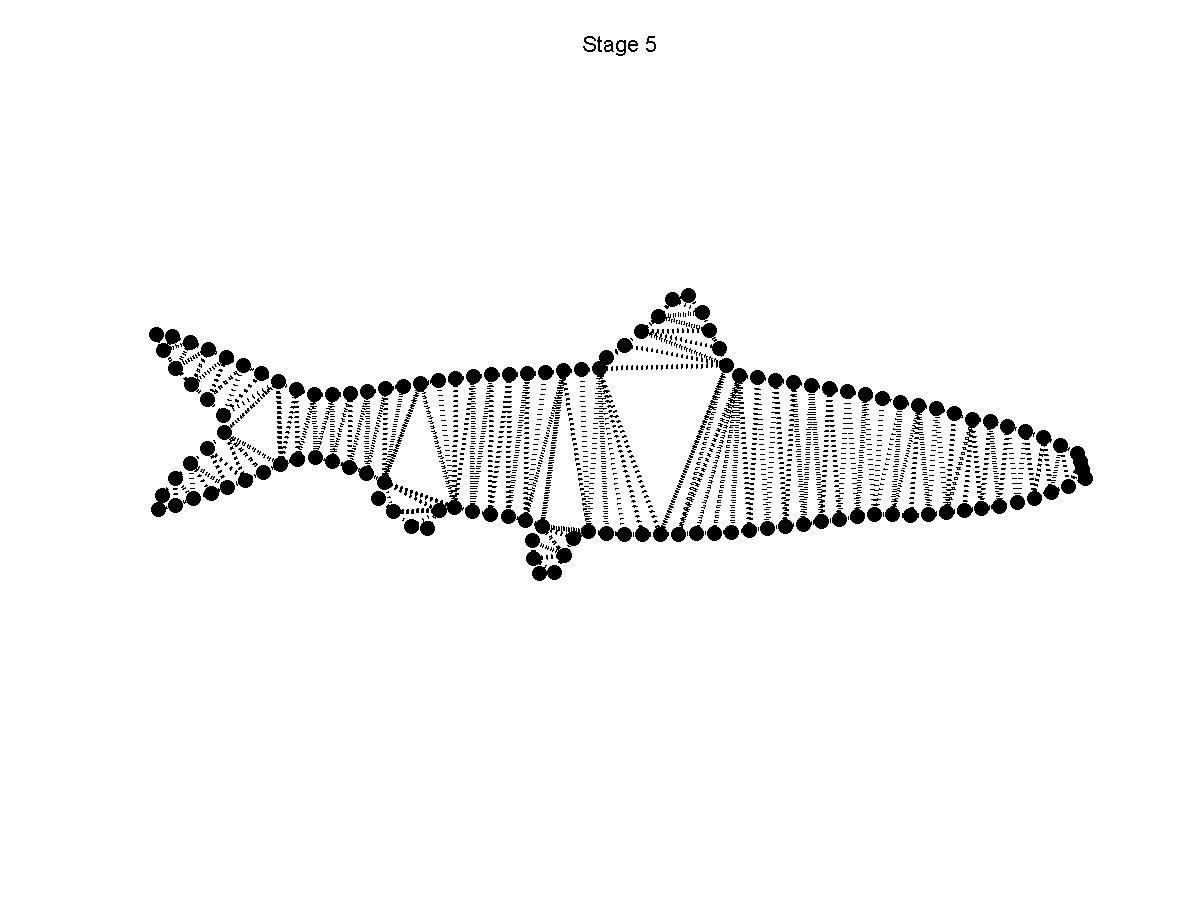}
\end{center}
\caption{Coarse to fine Delaunay triangulation used in the matching functional \eqref{eq:cr} for the shooting method. Triangulation is depicted on the subsets of $8,16,32,64$ points and the full set of 128 points of the fish.}
\label{fig:dt}
\end{figure}

We show the Delaunay triangulation for a fish shape at various refinement stages in Figure \ref{fig:dt}. Our choice of using a dyadic refinement strategy is not the only possibility, and we have chosen it mainly for convenience. In our results we employ between 4 and 5 refinement stages.

We summarize the full shooting procedure outlined in this section in Algorithm \ref{alg:method}.
\begin{algorithm}
\begin{algorithmic}
\STATE Input: number of teichons $N$ and $M$ ordered shape samples $z_m \in \C$
\STATE Determine number of refinement stages $S$: Section \ref{sec:algorithm-initial-guess}
\FOR{Stage $s=1,\ldots S$}
  \STATE Compute Delaunay triangulation and cross-ratios for $2^{(s-S)} M$ shape samples: Section \ref{sec:algorithm-delaunay} 
  \STATE Set initial guess for $p(0)$ as solution to previous stage: Section \ref{sec:algorithm-initial-guess}
  \WHILE{not converged}
    \STATE Compute objective $E_2$ \eqref{eq:cr} and the standard gradient \eqref{eq:E-chain-rule} using \eqref{sol_ode}, \eqref{eq:alpha-evolution}, \eqref{eq:cross-ratio}, and \eqref{eq:modified-ode}: Sections \ref{sec:algorithm-matching}, \ref{sec:algorithm-gradient} 
    \STATE Compute proper gradient using \eqref{eq:proper-gradient}: Sections \ref{sec:algorithm-gradient-projection}, \ref{sec:algorithm-natural-gradient}
    \STATE Perform either steepest descent \eqref{eq:p0-update}, or conjugate gradient descent \eqref{eq:cg-beta}, \eqref{eq:cg-update}: Section \ref{sec:algorithm-cg}
  \ENDWHILE
\ENDFOR
\STATE Output: $N$-teichon configuration $q(0)$, $p(0)$

\end{algorithmic}

\caption{Shooting algorithm for minimization of the cross-ratio objective \eqref{eq:cr}.}
\label{alg:method}
\end{algorithm}

\section{Numerical Results}
\label{s:numres}

In this section we present various numerical results that demonstrate the efficacy of our method. For all evolutions we employ $N=100$ teichons equispaced at $t=0$. Unless noted otherwise we iterate until the value of the objective \eqref{eq:cr} is no greater than $10^{-4}$, and frequently is $\mathcal{O}(10^{-6})$ at convergence. The number $M$ of landmarks $\alpha_m$ varies with the shape data, but usually takes the value $M=128$.

\subsection{Aspect ratio for an ellipse}
\begin{figure}
  \begin{center}
    \resizebox{1.0\textwidth}{!}{
      \includegraphics{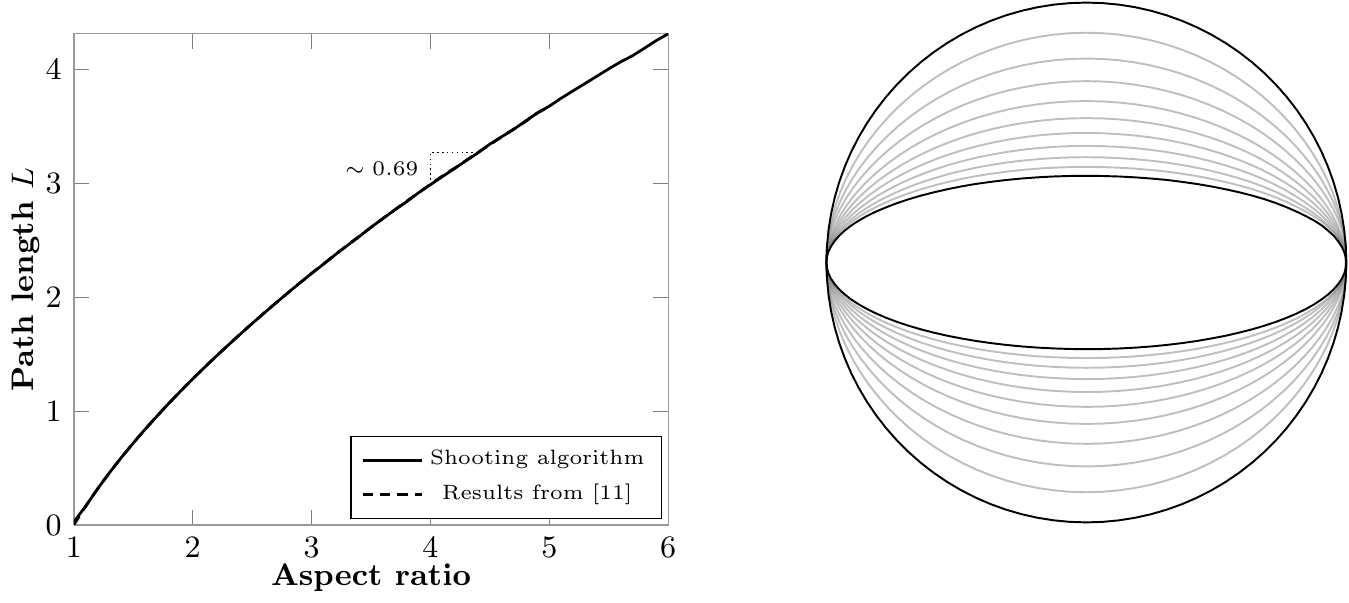}
    }
  \end{center}
  \caption{Left: $W P$ distance from circle to an ellipse vs. aspect ratio of the ellipse computed using the shooting algorithm presented in this paper (solid line) and using an energy minimization algorithm from \cite{feiszli_2012} (dashed line). (The lines overlap.) Right: Snapshots of the computed geodesic for aspect ratio 3 at equidistant points along the path.}
  \label{fig:ellipse-aspect-ratio}
\end{figure}

We first investigate the distance on $\uts$ between a circle and an ellipse of certain aspect ratio. With $M=100$ matching points, Figure \ref{fig:ellipse-aspect-ratio} graphs the distance on $T(1)$ between a circle and ellipses of aspect ratios from 1 to 6; these results are visually indistinguishable from those found with an energy minimization algorithm in \cite{feiszli_2012}, providing supporting evidence for the accuracy of the algorithm. The figure also suggests that the asymptotic ratio between the geodesic length and the aspect ratio of the ellipse is linear and the slope is approximately $0.69$. We note however that data for larger aspect ratios is necessary in order to verify this result. Unfortunately, numerical crowding prevents us from accurately computing geodesics for higher aspect ratios. 

%\tikzsetnextfilename{ellipse-aspect-ratio}

\subsection{Hyperbolicity test}
As has been mentioned in Section \ref{s:intro} all sectional curvatures of the Weil-Petersson metric are negative. In order to verify this numerically, we verify that the angle sum of a triangle on $T(1)$ with the $W P$ metric is less than $\pi$. We choose the three vertices of each triangle to be rotated ellipses of fixed aspect ratio. For a template ellipse with semimajor axis aligned with the horizontal axis, we choose the three vertices to be ellipses of rotations $0$, $2\pi/3$, and $4 \pi/3$; we label these vertices $s_1$, $s_2$, and $s_3$, respectively. We compute geodesics with $M=128$ matching points between two vertices of the triangle, and compute angles between these geodesics. 

\begin{figure}
  \begin{center}
    \resizebox{1.0\textwidth}{!}{
      \includegraphics{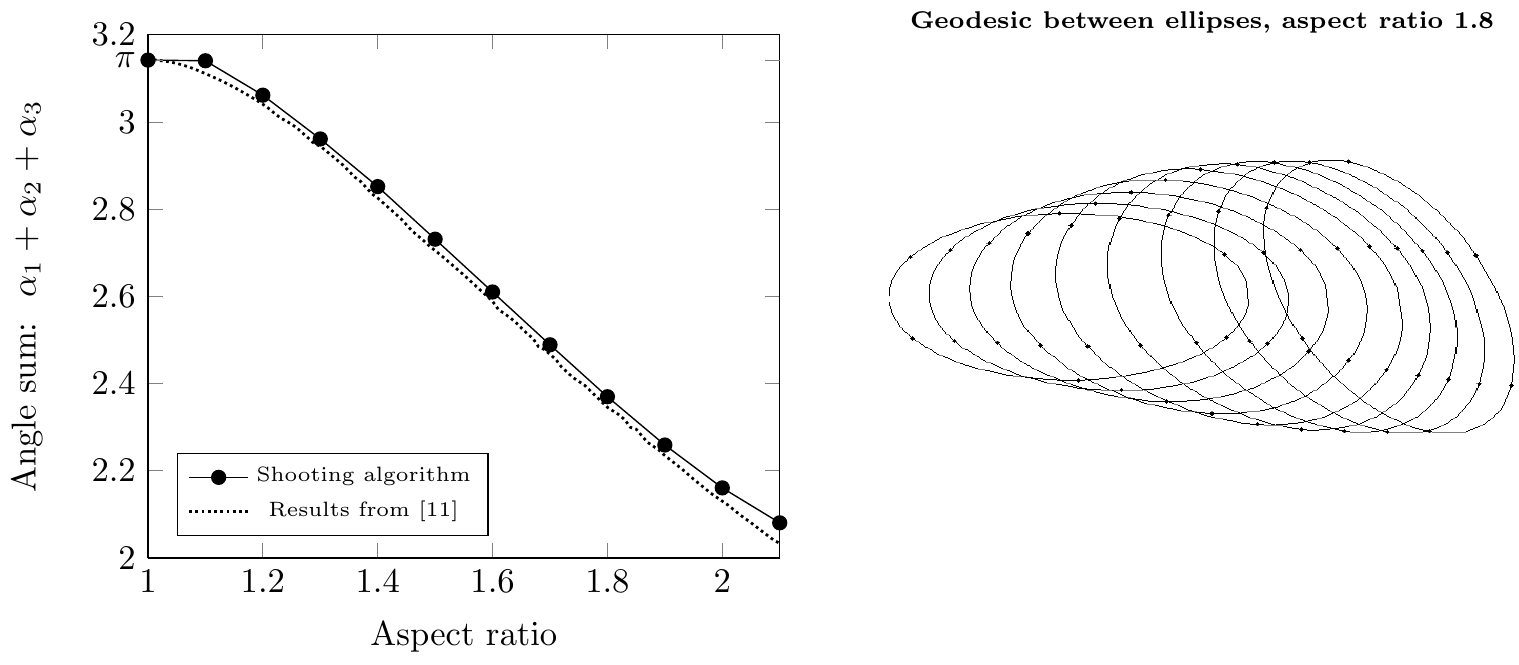}
    }
  \end{center}
  \caption{Left: Angle sums for triangles on $T(1)$ computed using the shooting algorithm of this paper (solid line) and the minimization algorithm from \cite{feiszli_2012} (dotted line). The three vertices for each triangle are formed by rotating a template ellipse of given aspect ratio. Right: sample evolution for ellipse of aspect ratio 1.6. An artificial shift is employed to make the evolution clearer.}
  \label{fig:ellipse-rotations}
\end{figure}

Let $v_{i,j} \in Hor$ denote the velocity field that pushes vertex $s_i$ to vertex $s_j$. The angles at the three vertices of this triangle are computed according to the formula:
%$$
%\alpha_1 = \frac{\langle v_{12},v_{13}\rangle_{WP}}{\|v_{12}\|_{WP} \|v_{13}\|_{WP}}.
%$$
\begin{align*}
  \alpha_i = \frac{\langle v_{i,i\oplus 1}, v_{i, i \ominus 1} \rangle_{W P}}{\|v_{i, i\oplus 1} \|_{W P} \|v_{i, i\ominus 1} \|_{W P}},
\end{align*}
where $\oplus$ and $\ominus$ denotes modular addition/subtraction on the set $\{1, 2, 3\}$. We then compute $\sum_{i=1}^3 \alpha_i$, which is the angle sum of the triangle. The graph of this sum versus the aspect ratio of template ellipse is depicted in Figure \ref{fig:ellipse-rotations}. The aspect ratio of ellipses varied from 1 to 2.2.

Once again we compare our results against those computed from an entirely different algorithm in \cite{feiszli_2012}. This comparison is shown in Figure \ref{fig:ellipse-rotations}, and in contrast to the previous test, there is now a noticeable difference in the results: The values differ by about $2\%$. We can attribute this difference to many factors. First, the sample points on the shape that we use in this algorithm are not the same points that are used in \cite{feiszli_2012}; this results in small differences in the conformal welds used for matching. Second, our algorithm uses $M=128$ matching points, whereas the energy minimization algorithm from \cite{feiszli_2012} used 150 points. Finally, our shooting method produces a numerically exact geodesic, but does not match the endpoint condition exactly; the energy minimization method from \cite{feiszli_2012} produces an approximate geodesic that matches the endpoint condition exactly to numerical precision. Thus, the small differences shown in the results are not altogether surprising.

%\tikzsetnextfilename{ellipse-rotations}

\subsection{A hippocampus slice}\label{ss:hippo}
\begin{figure}
  \begin{center}
    \resizebox{0.5\textwidth}{!}{
      \includegraphics{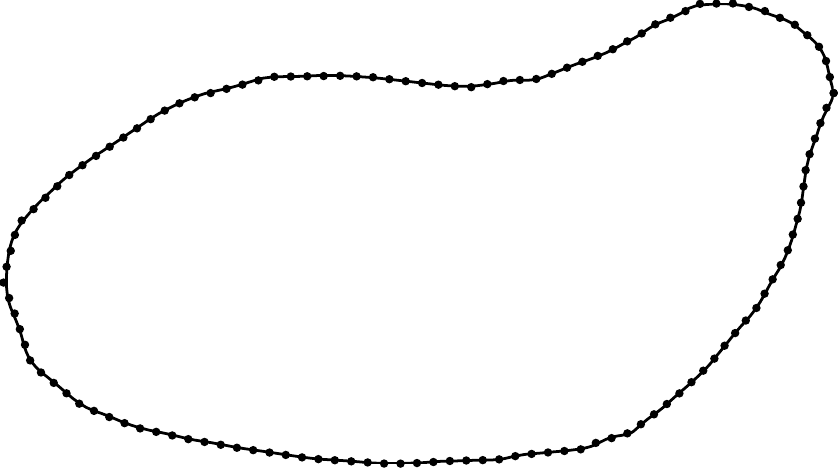}
    }
  \end{center}
  \caption{Terminal shape with shooting (solid line) versus landmarks used for the cross-ratio objective (solid dots).}
  \label{fig:hippo-target}
\end{figure}

\begin{figure}
  \begin{center}
    \resizebox{1.0\textwidth}{!}{
      \includegraphics{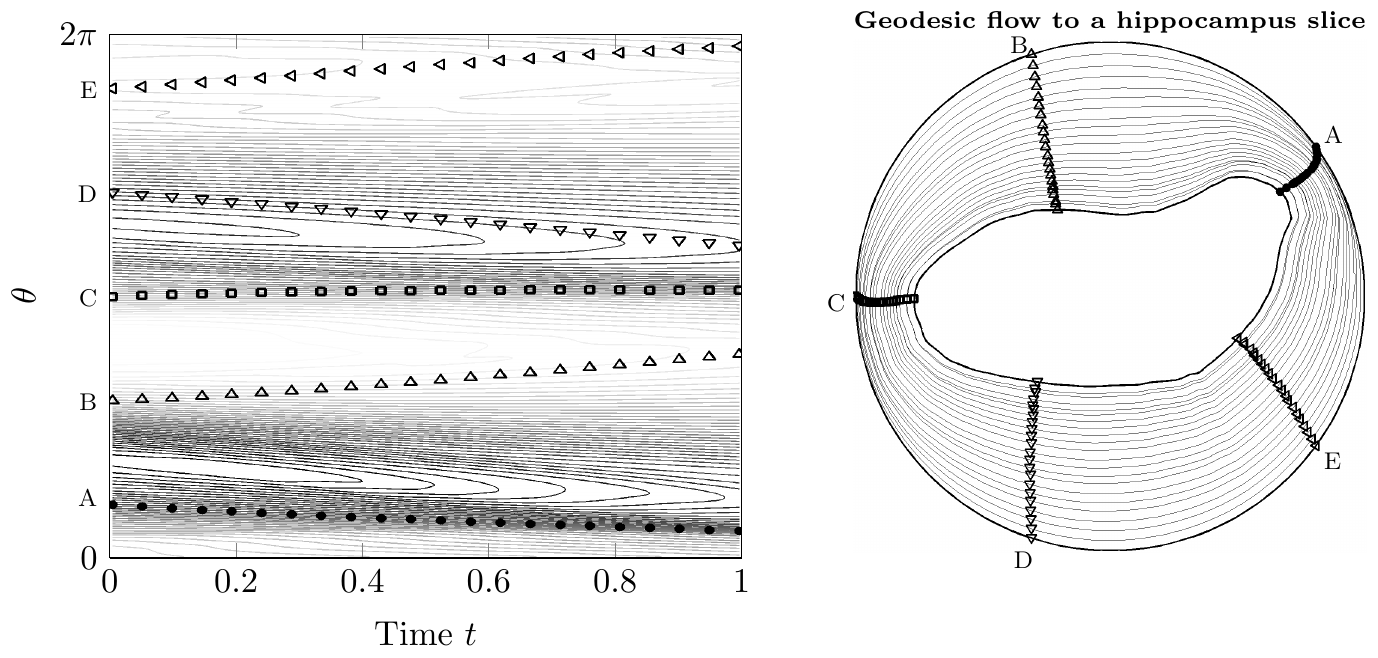}
    }
  \end{center}
  \caption{Left: Contour plot for velocity field $v(\theta,t)$ evolving a circle ($t=0$) to a hippocampus slice ($t=1$). Right: The resulting shape evolution. Shape snapshots are shown for $t=n/20$ for $n=0, \ldots, 20$. An artificial scaling is employed to make the evolution clearer. The path length on Teichm\"uller space is $1.968$.}
  \label{fig:hippo-evolving}
\end{figure}

We consider flowing from a circle to a planar slice of the human hippocampus in Figure \ref{fig:hippo-target} given by $M=128$ sample points. (For details, refer to \cite{kushnarev_2012}.) We shoot with 100 equidistant teichons and at algorithm termination the objective function is less than $10^{-4}$. Figure \ref{fig:hippo-target} shows that the shooting matches the target very well. Figure \ref{fig:hippo-evolving} displays the shape evolution along with a contour plot of the velocity field that pushes the circle to the hippocampus slice. We also track the evolution of 5 landmarks on the shape. The hippocampus slice is a relatively easy shape: the fingerprint does not exhibit crowding and so our computation of the geodesic flow (and the gradient) is accurate.

%\tikzsetnextfilename{hippo-target}

%\tikzsetnextfilename{hippo-evolving}

\subsection{Flow between shapes}\label{ss:bell}
\begin{figure}
  \begin{center}
    \resizebox{1.0\textwidth}{!}{
      \includegraphics{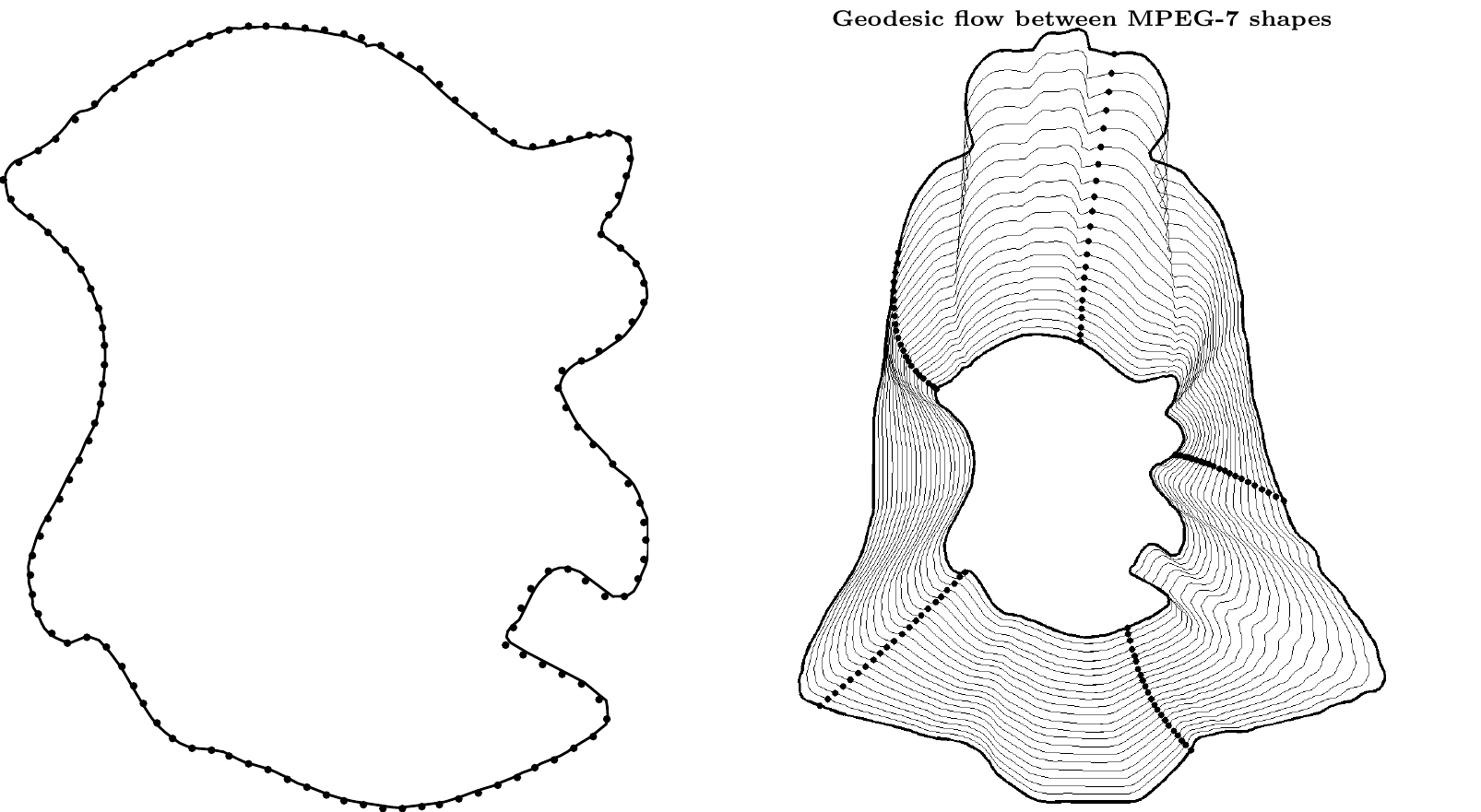}
    }
  \end{center}
  \caption{Shooting from a bell to a blob, two shapes in the MPEG-7 CE-Shape-1 database. Left: Terminal shape (blob) with shooting (solid line) versus landmarks used for the cross-ratio objective (solid dots). Right: Shape evolution with landmark flow (solid dots). An artificial scaling is employed to make the evolution clearer. The path length on Teichm\"uller space is 6.362.}
  \label{fig:bell-evolving}
\end{figure}

Our method does not rely on the initial shape being circular -- it is likewise possible to flow between non-circular shapes with no change to the algorithm. The MPEG-7 CE-shape-1 collection of planar shapes \cite{mpeg7} is a database of shapes commonly used in classification routines. Our immediate goal here is not classification, but to illustrate the applicability of our algorithm to realistic shapes, we choose two shapes with non-crowded welding maps from this database and show the Teichm\"uller evolution between them obtained from the shooting algorithm with $M=128$ matching points, see Figure \ref{fig:bell-evolving}.
%\tikzsetnextfilename{bell-evolving}

\subsection{A fish contour}\label{ss:other-shapes}
\begin{figure}
  \begin{center}
    \resizebox{1.0\textwidth}{!}{
      \includegraphics{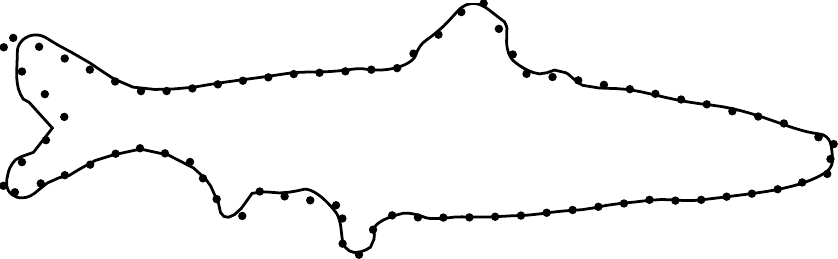}
    }
  \end{center}
  \caption{Terminal shape with shooting (solid line) versus landmarks used for the cross-ratio objective (solid dots).}
  \label{fig:fish-target}
\end{figure}

\begin{figure}
  \begin{center}
    \resizebox{1.0\textwidth}{!}{
      \includegraphics{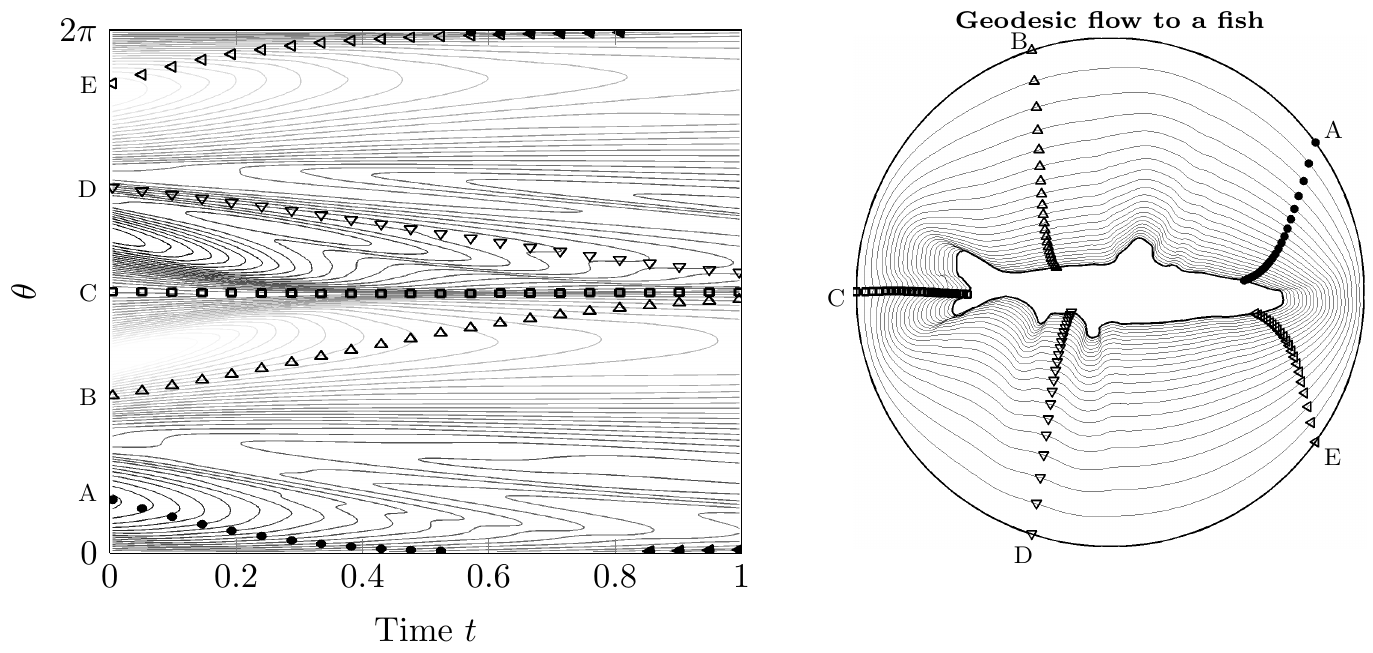}
    }
  \end{center}
  \caption{Left: Contour plot for velocity field $v(\theta,t)$ evolving a circle ($t=0$) to a fish ($t=1$). Right: The resulting shape evolution. Shape snapshots are shown for $t=n/20$ for $n=0, \ldots, 20$. An artificial scaling is employed to make the evolution clearer. The path length on Teichm\"uller space is $7.152$.}
  \label{fig:fish-evolving}
\end{figure}

We finally consider a more complex shape: an outline of the fish given in Figure \ref{fig:fish-target} from $M=128$ sample points. The initial shape is a circle, and the initial teichon configuration is given by 100 equidistant teichons. The result of the matching is given in Figure \ref{fig:fish-target}. As with the hippocampus slices, we show the resulting $W P$ geodesic from a circle to the fish on the right panel of Figure \ref{fig:fish-evolving}, and a contour plot of the velocity is shown on the left panel. The objective value at termination of the algorithm is $5 \times 10^{-4}$. We see immediately that there are limitations to this algorithm: the welding map for the fish outline suffers from severe crowding. 

For the particular chart we have chosen for the fingerprint $\psi_1$, the landmarks on the tail of the fish are separated by a distance of $\mathcal{O}(10^{-8})$. This crowding of the shape landmarks $\alpha_m$ implies a similar crowding of the teichon positions $q_k$ that push the landmarks. When teichon positions $q_k$ and $q_{k+1}$ are $\mathcal{O}(10^{-8})$, we can no longer accurately integrate EPDiff. To understand why, consider first the Gram matrix $G$ for teichon positions $q_k$ (see \eqref{eq:gradient-projection}), which has entries $G_{k,l} = G(q_k - q_l)$, with $G(\cdot)$ being the $W P$ Green's function \eqref{eq:greens-function}. From the explicit form of the Green's function, one can show that for small arguments $\Delta \theta$,
\begin{align*}
  G(\Delta \theta) \approx \frac{1}{2} + \Delta \theta^2 \left[ \log(\Delta \theta^2) - \frac{3}{4}\right].
\end{align*}
Then when $q_k - q_{k+1}$ is $\mathcal{O}(10^{-8})$, we have $G_{k,k+1} = G(q_k - q_{k+1}) \approx \frac{1}{2} + \mathcal{O}(10^{-16})$, and in double-precision arithmetic where we have implemented this code, floating-point truncation error causes this matrix entry to coincide with $G_{k,k} = G(0) = \frac{1}{2}$. This means that the Gram matrix is singular to numerical precision. Thus we cannot accurately evaluate the right-hand side of EPDiff given by \eqref{sol_ode} and \eqref{eq:alpha-evolution}, and also cannot accurately compute the gradient using the EPDiff-derived system \eqref{eq:modified-ode}. Therefore the algorithm begins to break down at this point in the sense that we cannot resolve features that require teichons to flow so close to one another on $S^1$.

In general, when teichons flow very close to one another (relative to machine precision), we observe signature failures of the algorithm due to finite precision through various diagnostics:
\begin{itemize}
  \item the computed $W P$ norm of the $N$-teichon is not constant in time $t$, or becomes negative,
  \item teichon locations $q_k$ (or landmarks $\alpha_m$) cross each other,
  \item stepping in the direction of the proper gradient does not decrease the objective.
\end{itemize}
The first two issues can normally be ameliorated by decreasing the $\Delta t$ time-stepping parameter used to integrate \eqref{sol_ode}, but the third issue is usually difficult to resolve in an automated fashion.

%This value of $10^{-8}$ is a critical point for the algorithm: the Green's function \eqref{eq:greens-function} for the $W P$ operator has a quadratic maximum at $\theta=0$ -- this implies that the Gram matrix is numerically singular for two points $\theta_1$ and $\theta_2$ whose separation $|\theta_1 - \theta_2|$ is smaller than the square root of machine precision $\sqrt{\epsilon_{\mathrm{mach}}}$. Since we have implemented this method in double-precision arithmetic, $\sqrt{\epsilon_{\mathrm{mach}}}$ works out to be about $10^{-8}$. Thus, we are unable to precisely match the fish welding map because numerical precision prevents us from properly resolving landmarks (and thus features) whose separation is $\mathcal{O}(10^{-8})$.

%\tikzsetnextfilename{fish-target}

%\tikzsetnextfilename{fish-evolving}

\section{Conclusions}
In this paper we have demonstrated an efficient method for computing geodesics on the coset space $\PSL_2(\R)\bs\Diff(S^1)$, a dense subset of the universal Teichm\"uller space $T(1)$, with the Weil-Petersson metric via shooting. The geodesics are found by approximating the velocity field with an ansatz solution, an $N$-Teichon. The fact that an $N$-Teichon solution remains an $N$-Teichon under geodesic flow allows us to accurately compute these geodesics. A matching term is employed to guide the initial guess; cross-ratios make up the matching term to correctly identify disparities between equivalence classes on the coset space $\PSL_2(\R)\bs\Diff(S^1)$. We are able to use a nonlinear optimization algorithm to converge to geodesics on this space. However, our method still suffers from the well-known crowding phenomenon, which prevents us from computing geodesics between shapes with elongated features. Even if a crowded welding map can be accurately computed, the geodesic evolution becomes inaccurate when particles flow to within $\sqrt{\epsilon_{\mathrm{mach}}}$, where $\epsilon_{\mathrm{mach}}$ is machine precision for floating-point computations. Nevertheless, there is a wide range of non-crowded shapes for which our algorithm is effective.

To our knowledge this is one of only a few numerical algorithms that can reliably compute the Weil-Petersson geodesics on the $T(1)$. It performs much better than the energy minimization method proposed in \cite{sharon_2d-shape_2006}, and is competitive with the recent approach \cite{feiszli_2012}. An additional advantage to the proposed method is that it is a shooting method: a proper geodesic is always produced (up to the precision of the forward integration scheme). Future work will employ this method for consistent comparison of shapes in a database. The uniqueness of geodesics implies that consistency is ensured by the unique initial momentum that is assigned to each shape via the tangent space linearization. Moreover, this allows us to find unique Karcher mean and perform well-posed statistics on the shape space \cite{karcher}. In \cite{kushnarev_2012} we have employed the described method to study the database of hippocampus of patients with dementia and healthy controls.

\section*{Acknowledgments}
The authors would like to thank Prof.~David Mumford, Prof.~Darryl Holm and Asst Prof.~Anqi Qiu for their insightful suggestions, discussions and comments.

\bibliography{shooting_teichon}
\bibliographystyle{abbrv}

\end{document}